\documentclass[12pt]{iopart}
\usepackage{iopams}
\usepackage{graphicx}
\usepackage{clrscode}
\usepackage{subfigure}
\usepackage{color}
\graphicspath{{figs/}}

\def\(({\left(}
\def\)){\right)}
\def\[[{\left[}
\def\]]{\right]}

\newcommand{\atanh}{{\rm{atanh}}}
\newcommand{\be}{\begin{equation}}
\newcommand{\ee}{\end{equation}}
\newcommand{\bea}{\begin{eqnarray}}
\newcommand{\eea}{\end{eqnarray}}

\newcommand{\s}{x}

\newcommand{\apriori}{{\it a priori }}
\newcommand{\aposteriori}{{\it a posteriori }}
\newcommand{\vecs}{{\mathbf{\s}}}
\newcommand{\vecy}{{\mathbf{y}}}

\begin{document}

\title{Belief Propagation Reconstruction for Discrete Tomography}

\author{E. Gouillart$^1$, F. Krzakala$^2$, M. M\'ezard$^3$ and L. Zdeborov\'a$^4$}

\address{$^1$Surface du Verre et Interfaces, UMR 125
  CNRS/Saint-Gobain, 93303 Aubervilliers, France\\ $^2$ CNRS and ESPCI
  ParisTech, 10 rue Vauquelin, UMR 7083 Gulliver, Paris 75005,
  France\\$^3$ Ecole normale sup\'erieure, 45 rue d'Ulm, 75005 Paris, and LPTMS, Univ. Paris-Sud/CNRS, B\^at. 100, 91405
  Orsay, France\\$^4$ Institut de Physique Th\'eorique, IPhT, CEA
  Saclay, and URA 2306, CNRS, 91191 Gif-sur-Yvette, France}
\ead{emmanuelle.gouillart@nsup.org}

\begin{abstract}

  We consider the reconstruction of a two-dimensional discrete image from
a set of tomographic measurements corresponding to the Radon projection.
Assuming that the image has a structure where neighbouring pixels have a
larger probability to take the same value, we follow a Bayesian
approach and introduce a fast message-passing reconstruction algorithm
based on belief propagation. For numerical results, we specialize to the
case of binary tomography. We test the algorithm on binary synthetic
images with different length scales and compare our results against a
more usual convex optimization approach. We investigate the
reconstruction error as a function of the number of tomographic
measurements, corresponding to the number of projection angles. The
belief propagation algorithm turns out to be more efficient than the
convex-optimization algorithm, both in terms of recovery bounds for
noise-free projections, and in terms of reconstruction quality when
moderate Gaussian noise is added to the projections.

\end{abstract}

\submitto{\IP}

\maketitle

\section{Introduction}

\begin{figure}
\begin{center}
\centerline{\includegraphics[width=0.7\textwidth]{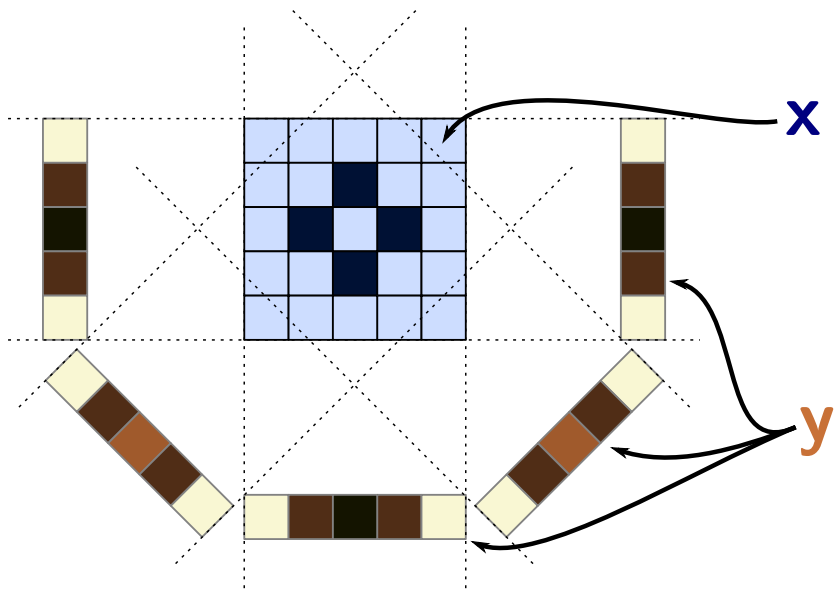}}
\end{center}
\caption{Geometry of the tomography problem: the tomographic measurements
$\mathbf{y}$ are line sums of the image $\mathbf{x}$ along different
angles.
\label{fig:tomo_diagram}}
\end{figure}

X-ray computed tomography (CT)~\cite{Slaney1988, Herman2009} is a
classical 3-D imaging technique in materials science or for medical
applications. The X-ray beam transmitted through the sample is recorded
on a planar detector, for different angles of incidence of the beam on
the sample. Using the Beer-Lambert law for photons
absorption, the reconstruction of the absorption $\mathbf{x}$ of the
sample from the measurements $\mathbf{y}$ is a linear problem

\be
\mathbf{y} = \mathbf{F}\mathbf{x} + \mathbf{w},
\label{eq:tomo_linear}
\ee
with $\mathbf{F}$ the tomography projection matrix and $\mathbf{w}$ the
measure noise. For a parallel X-ray beam, each pixel line of the detector
measures the total absorption integrated along a horizontal light-ray
through the sample. The geometry of tomographic measurements in shown on
the schematic of Fig.~\ref{fig:tomo_diagram}. For a detector line of $L$
pixels, a square image of $N=L^2$ pixels has to be reconstructed from
$M=L\times n_{\theta}$ measurements, where $n_\theta$ is the number of
angles, also called number of projections.

Classical reconstruction algorithms such as filtered back-projection
(FBP)~\cite{Shepp1974,Ramachandran1971,Slaney1988} require the linear
system (\ref{eq:tomo_linear}) to be sufficiently determined. Therefore,
the number of angles $n_\theta$ needed is generally close to the number
of pixels along a line, $L$. Nevertheless, many applications would
benefit from being able to reconstruct from a smaller number of angles.
In medical imaging for example, reducing the total X-ray dose absorbed by
the patient is highly desirable. For in-situ tomography where the sample
is evolving while it is imaged~\cite{Baruchel2006,Buffiere2010}, reducing
the number of angles increases the acquisition speed, hence the
time-resolution of the observations. 

Filtered back-projection fares poorly for a limited number of angles, and
algorithms incorporating \emph{prior information} on the sample are
needed to improve the quality of the reconstruction. For example, a
positivity constraint on the absorption incorporated in the Algebraic
Reconstruction Technique (ART) improves the result of the
reconstruction~\cite{Gordon1970}. More specific priors can also be used.
In particular, samples with a discrete number of phases and absorption
values are frequently encountered in materials science. Using this
information can make up for the lack of measurements and result in a
satisfying reconstruction. The field of \emph{discrete
tomography}~\cite{Discrete09} encompasses the large class of algorithms
that have been proposed to solve this problem. Various approaches have
been proposed, from heuristic methods~\cite{Batenburg2011} to the
approximate resolution of the combinatorial problem~\cite{OldBayes}, or
of its convex relaxation~\cite{Sidky2008}. 

In this article, we propose a new algorithm for discrete tomography, that
estimates the marginal probability for every pixel value. We use a
Bayesian approach in order to write the \emph{a-posteriori} probability
distribution of discrete images, and a fast approximate message passing
scheme relying on belief propagation to find the marginal probabilities.
Within this framework, constraints imposed by the measurements are
satisfied by iteratively solving a graphical model along the light-rays
crossing the sample, known as the Potts model in statistical physics. For
the specific case of binary images, this amounts to solving 1-D Ising
models for all light-rays and measurements. The Potts and Ising models
are specific examples of Markov random fields (MRF). In the case of image
segmentation, these models can be solved using graph-cuts
methods~\cite{Boykov2001}.

We start this article with a short review of the state of the art on
discrete tomography, where we introduce the convex algorithm based on
total variation, used for comparison with our algorithm. Then we move on
to the presentation of our belief-propagation tomographic reconstruction
algorithm (BP-tomo) and its relation to the field of statistical physics.
Numerical results on synthetic data are presented for the case of binary
images. We establish empirically the recovery phase diagram by decreasing
the number of projections and determining the undersampling rate
separating successful and failed reconstruction; it is found that the
transition corresponds to an undersampling rate of the order of the
density of gradients in the image. This bound on undersampling rates is a
significant improvement over performances of traditional compressed
sensing. In the realistic case of noisy measurements, reconstructions
obtained with BP-tomo are compared with results obtained with a
constrained total-variation minimization. For moderate noise, BP-tomo is
found to give a much more accurate reconstruction above the critical
undersampling rate of the noise-free case.

\newpage
\section{State of the art in discrete tomography \label{sec:state_of_art}}

Discrete tomography algorithms perform at the same time the tomographic
\emph{reconstruction} and the \emph{segmentation} of the reconstructed
image into objects with a finite number of known absorption values. In
addition to the constraints on the discrete pixel levels, that strongly
restrict the space of possible images, additional constraints on the
gradients of the image are often enforced in order to obtain smooth
images. A large range of methods has been proposed for discrete
tomography, for a review see~\cite{Discrete09}. Here, we only cite a few
representative examples, in order to discuss briefly their scope of
application and performance.

Theoretical studies~\cite{Gardner1995, Discrete09} have addressed the
problem of the minimal number of angles needed to reconstruct a binary
image, in the noise-free case. An important result of \emph{geometric
tomography}~\cite{Gardner1995} is that four different projections are
enough to reconstruct uniquely a uniform convex
object~\cite{Gardner1980}. For more complicated images, degenerate binary
solutions may exist for the same set of projections, but it is possible
to bound the distance between two solutions using the projection error of
the thresholded filtered-back projection solution~\cite{Batenburg2011b}.
It is also possible~\cite{Aharoni1997} to determine the set of pixels
uniquely determined (in the absence of noise) by their projections. 

As for reconstruction algorithms, a first class of
methods~\cite{OldBayes,Liao2004,Liao2005,Liao2006,Djafari2008} adopts a Bayesian
approach in order to make the most of the available information on the
sample. The \emph{a posteriori} probability distribution of a discrete
image $\mathbf{{ {x}}}$ given measurements $\mathbf{y}$ is
\be 
P(\mathbf{x}|\mathbf{y})=\frac{P(\mathbf{y}|\mathbf{x})
P_0(\mathbf{x})}{P(\mathbf{y})},
\ee
with $P_0(\mathbf{x})$ the \emph{a priori} probability distribution
on the space of images. In the Maximum A Posteriori (MAP) setting, the
image maximizing the posterior distribution is searched for:
\be
x_{\mathrm{MAP}} = \mathrm{argmax}\quad P(\mathbf{x}|\mathbf{y}).
\ee
Finding $x_{\mathrm{MAP}}$ is an NP-hard combinatorial problem. An
MCMC (Monte-Carlo Markov Chain) algorithm is used
in~\cite{OldBayes,Liao2004,Liao2005,Liao2006} in order to sample the
Gibbs distribution and to obtain an approximation of
$x_{\mathrm{MAP}}$. Gibbs sampling being very costly in terms of
numbers of iterations, such algorithms are limited to small images.
Therefore, these methods are more suitable for techniques with lower
spatial resolution than computed tomography, such as positron emission
tomography (PET) for example. For faster computations, a variational
Bayes approach was used in~\cite{Djafari2008}. The Bayesian approach is very
flexible in terms of prior distribution on the images, allowing for a
tuned weighting of local binary patterns in the image~\cite{Liao2004}.
For binary images, a popular model for smooth images is the Ising model
from statistical physics, that we shall use as well in our
belief-propagation algorithm.

A second class of
algorithms~\cite{Delaney1998,Weber2003,Schule2005,Weber2006,Weber2006b}
consider relaxations of the MAP estimation falling into the scope
of \emph{linear optimization}. In this setting, the searched-for image
$\vecs$ is found by solving a set of linear constraints, such as the
constraint $\mathbf{y} = \mathbf{F}\mathbf{x}$, or an interval
constraint. For mixtures of convex and concave constraints, methods
from DC (Difference of Convex functions) programming are used to converge
to a local minimum of the objective
function~\cite{Schule2005,Weber2005,Weber2006,Weber2006b}.

Since the seminal papers by Candes, Romberg, Tao and
Donoho~\cite{Candes2006, Candes2006b, Donoho2006}, the field of
\emph{compressed sensing} has breathed a new lease of life into the use
of convex optimization for finding sparse solutions to underdetermined
inverse problems. Compressed sensing proposes to find sparse solutions
using a regularization with the convex $\ell_1$ norm, instead of the
$\ell_0$ norm measuring sparsity: problems can then be solved using
classical techniques of convex optimization~\cite{Boyd2004,
Combettes2011}. In general, convex optimization methods can be
computationally less intensive than other types of methods, hence more
suitable for large images. Candes et al.~\cite{Candes2006b} proved under
certain hypotheses that perfect reconstruction can be obtained with the
$\ell_1$ norm in the noiseless case, and that the error is at most
proportional to the noise level. These solid mathematical foundations
have contributed greatly to the success of compressed sensing. Although
tomographic projection does not satisfy the restrictive hypotheses of the
compressed sensing setting, compressed-sensing techniques have proven to
work well for tomographic reconstruction. In fact, the founding paper of
compressed sensing~\cite{Candes2006} demonstrated perfect reconstruction
with $\ell_1$ regularization using the example of incomplete tomographic
projections of the piecewise-constant Shepp-Logan
phantom~\cite{Shepp1974}. For the piecewise-constant images of discrete
tomography, the regularization is on the $\ell_1$ norm of the gradient of
the image, that is its total variation (TV)
\be
\mathrm{TV}(\mathbf{x}) = \sum_{i=0}^N
|\nabla \mathbf{x}_i|\ ,
\ee
where $\nabla \mathbf{x}_i$ is the discrete gradient of the image
evaluated at pixel $i$.
Therefore, several papers~\cite{Song2007, Sidky2008, Herman2008,
Tang2009, Sidky2010, Jia2010, Anthoine2011} have used the following
convex minimization to reconstruct discrete tomography images:
\be
\mathbf{x}_{\mathrm{TV}} = \mathrm{argmin}_{\mathbf{x}}\quad \frac{1}{2}
|| \mathbf{y} - \mathbf{F}\mathbf{x} ||^2 + \beta\mathrm{TV}(\mathbf{x})
\label{eq:tv_reg}
\ee 
where $|| \mathbf{y} - \mathbf{F}\mathbf{x} ||^2$ is the $\ell_2$ data
fidelity term (corresponding to a Gaussian noise prior in the MAP
setting), and the parameter $\beta$ controls the amount of
regularization. It has been shown~\cite{Song2007, Sidky2008, Herman2008,
Tang2009, Sidky2010, Jia2010, Anthoine2011} that the total variation
reconstruction greatly improves the quality of the reconstruction of
discrete images for a reduced number of projections. There is often a
trade-off (controlled by the parameter $\beta$) between removing the
noise and smoothing out small features~\cite{Tang2009}. For enforcing
spatial regularization, it is also possible to use frame
representations~\cite{Pustelnik2009, Pustelnik2010} -- especially for
medical images that are often more complex than piecewise-constant -- or
to combine frames and total variation~\cite{Pustelnik2010}. 

In this paper, we use a total-variation minimization algorithm to benchmark the
quality of the reconstruction of the belief-propagation algorithm.
Compared to Eq. (\ref{eq:tv_reg}), we also constrain the solution to have
values between the minimal discrete value $a$ and the maximal value
$b$~\cite{Pustelnik2010}.
We solve
\be
\mathbf{x}_{\mathrm{cTV}} = \mathrm{argmin}_{\mathbf{x}}\quad \frac{1}{2}
|| \mathbf{y} - \mathbf{F}\mathbf{x} ||^2 + \beta\mathrm{TV}(\mathbf{x})
+ \mathcal{I}_{[a,b]}(\mathbf{x}), 
\label{eq:tv_interval}
\ee
where $\mathcal{I}_{[a,b]}(\mathbf{x})$ is the multivariate indicator
function of the interval $[a,b]$. In order to solve the minimization
problem, we use proximal methods~\cite{Combettes2011} that perform
subgradient iterations on non-differentiable terms. Several
methods~\cite{Pustelnik2009, Raguet2011} are suited for a sum of
non-differentiable functionals (such as $\mathcal{I}_{[a,b]}$ and
$\mathrm{TV}$), we use here the generalized forward-backward
splitting~\cite{Raguet2011}. We have verified experimentally that with
the additional interval constraint, better reconstruction results are
obtained with Eq. (\ref{eq:tv_interval}) than with Eq. (\ref{eq:tv_reg}).
For the total variation proximal operator, we use the second-order FISTA
scheme of Beck and Teboulle~\cite{Beck2009} for the isotropic total
variation.
 
Finally, several heuristic algorithms have been proposed for discrete
tomography. One of the simplest~\cite{Myers2008} consists in alternating
gradient iterations on the data fidelity term and clipping (thresholding)
continuous values on the discrete labels.
Batenburg~\cite{Batenburg2007} proposed a network-flow algorithm where
constraints imposed by pairs of directions are satisfied iteratively for
binary pixel values. For binary image, Batenburg also introduced
DART~\cite{Batenburg2011}, a variation on the gradient descent - clipping
method: pixels are progressively labeled from the interior of objects to
the more uncertain boundaries. This algorithm takes into account the
spatial continuity of the objects as well as their binary nature. Along
the same lines, another algorithm was recently proposed by Roux et
al.~\cite{Roux2012}, where the reconstruction is done for the continuous
probability of pixels values, reducing the risk of clipping too early the
value of a pixel. 

Compared to the existing literature on discrete tomography, we shall see
that our belief propagation algorithm combines the full power of the
Bayesian approach and the efficiency of iterative reconstruction
algorithms. 

\newpage
\section{Image reconstruction as a statistical physics problem 
\label{sec:statphys}}

\subsection{The Bayesian approach}
Let us first fix the notations. The image that we want to reconstruct sits on a
two-dimensional grid, and each of its $N$ pixels $i$ takes a value in
a $q$-ary alphabet: $x_i=\chi_1,\ldots,\chi_q$. For instance, binary
tomography corresponds to the case where $q=2$, and $\chi_1=-1$, $\chi_2=1$.  We suppose that $M$ tomographic measurements, $\mu=1
\ldots M$, have been made. Each such measurement gives a number
$y_{\mu}$ which is the sum of variables along a light ray $\mu$:
\be y_{\mu}=\sum_{i \in {\mu}} \s_i + w_\mu \ee 
where ${i \in {\mu}}$ denotes the set of all variables $\s_i$ that
belong to the measurement line $\mu$, and $w_\mu$ is the additive noise on measure
$\mu$. Notice that, in a slightly more
general setting, one would measure $y_{\mu}=\sum_{i \in {\mu}} F_{\mu
  i}\s_i + w_\mu$, where $F_{\mu i}$ is a geometrical factor giving the
fraction of pixel $i$ that is covered by the ray corresponding to
measurement $\mu$~\cite{Joseph1982}. This more general setting can be handled as well by
our method. Furthermore, it does not make much difference when the
image has a structure in domains, as soon as the typical domain size is
much larger than one pixel. To keep notations simple we thus keep
hereafter to the case where $F_{\mu i}=1$, the straightforward extension to
more general cases is left to the reader.

The problem is now to
reconstruct the original image (that is the N-dimensional vector
$\vecs$) from the values of the projection (the M-dimensional vector
$\vecy$).

Here, we shall adopt a Bayesian approach. As a prior information, we use
the fact that the image to measure is not a random one, but has a
structure in domains, such that if a pixel $i$ has the value $\chi_r$,
then the neighboring pixel has a large probability to have the same
value. We shall enforce this bias in each direction $\mu$ in which a
measurement has been made. This is a very natural prior for a large class
of images. Following the Bayesian approach, our goal is thus to sample
from the \aposteriori distribution given by
\bea P(\vecs|\vecy)=\frac{P(\vecy|\vecs)  P_0(\vecs)}{P(\vecy)} 
= \frac{1}{Z(\vecy)} {\prod_{\mu}\delta\(({y_{\mu} - \sum_{i \in {\mu}}
    \s_i}\))} P_0(\vecs) \label{Bayes_Intro} \eea
where, following the convention in statistical physics, we have
rewritten the normalization factor of $ P(\vecs|\vecy)$ as $Z(\vecy)$. (In
the following, we shall incorporate all normalization factors into $Z$.)
Here the $\delta$ function implements the equality $y_{\mu} = \sum_{i
  \in {\mu}} \s_i$ if the measurement is perfect ($w_\mu=0$). In the case of noisy
measurement, this $\delta$ function should be smoothed in order to
take into account our information on the noise. For instance, when the
measurement has Gaussian noise of variance $\Delta$ one could use a
function $\delta(x)=\exp(-x^2/(2\Delta))/\sqrt{2\pi \Delta}$.

Let us denote by $p_{\mu}$ the \apriori probability for neighboring
pixels to have the same value along the line $\mu$. Then
\bea 
P_0(\vecs) &=& \prod_{\mu} \prod_{(ij) \in \mu } p_{\mu}^{\delta_{\s_i,\s_j}}
(1-p_{\mu})^{1-\delta_{\s_i,\s_j}} = \prod_{\mu} (1-p_\mu)^{N_\mu}  \prod_{\mu}
e^{\sum_{ (ij)\in \mu} \log\left(\frac{ p_{\mu} }{1-p_{\mu}}
  \right)\delta_{\s_i,\s_j}} \nonumber \\ &=& {\rm{cst}}   \prod_{\mu}
e^{\sum_{(ij) \in \mu}  J_{\mu}\delta_{\s_i,\s_j}}\label{P0def}
\eea
where $N_\mu$ is the number of pixels along line $\mu$, $(ij)\in \mu$ denote pairs of neighboring pixels along direction
$\mu$ (for instance, in the case of line $\mu=3$ in Fig.\ref{Howto},
there would be an interaction between pixels (i) and (ii) and another
one between (ii) and (iii)) and $\rm{cst}$ is a constant that does not
depend on the values of the variables, and that we shall thus incorporate
into the normalization factor through a redefinition of $Z(\vecy)$.  In
(\ref{P0def}) we have introduced the so-called ``coupling constant''
$J_{\mu} = \log{\frac{ p_{\mu} }{1-p_{\mu}}}$.

With these notations, the \aposteriori probability can thus be written
as a product of $M$ terms, each corresponding to one tomographic
projection:
\bea
P(\vecs|\vecy) = \frac 1Z \prod_{\mu=1} ^M\left[
\delta\((y_{\mu} - \sum_{i \in  {\mu}}    \s_i\)) e^{J_{\mu} \sum_{(ij)
    \in \mu}  \delta_{\s_i,\s_j}} \right]
\label{Bayes}
\eea

The value of $J_{\mu}$ can be fixed based on previous knowledge of the
spatial structure of the sample. It would also be possible to learn its
value in the course of the algorithm, by maximizing the
probability of $J_{\mu}$ given the measures $\mathbf{y}$ (this amounts to
maximizing the partition function $Z$ in Eq. (\ref{Bayes})). This can be
done for instance using the Expectation Maximization
approach~\cite{Dempster}.

\subsection{A message passing approach: Belief Propagation}
Sampling from the distribution in Eq. (\ref{Bayes_Intro}) is a
notoriously difficult problem \cite{OldBayes,Discrete09}. We shall
here use an approach based on a message-passing algorithm called
``Belief-Propagation'' (BP)
\cite{Pearl82,KschischangFrey01,YedidiaFreeman03,MezardMontanari09}. It
is not an exact algorithm (at least for the present problem), in the
sense that it is not guaranteed to sample correctly from the
distribution (\ref{Bayes_Intro}). However, we shall see that the
approximate sampling obtained from BP performs empirically very well.

\begin{figure}
\begin{center}
  \includegraphics[width=14cm]{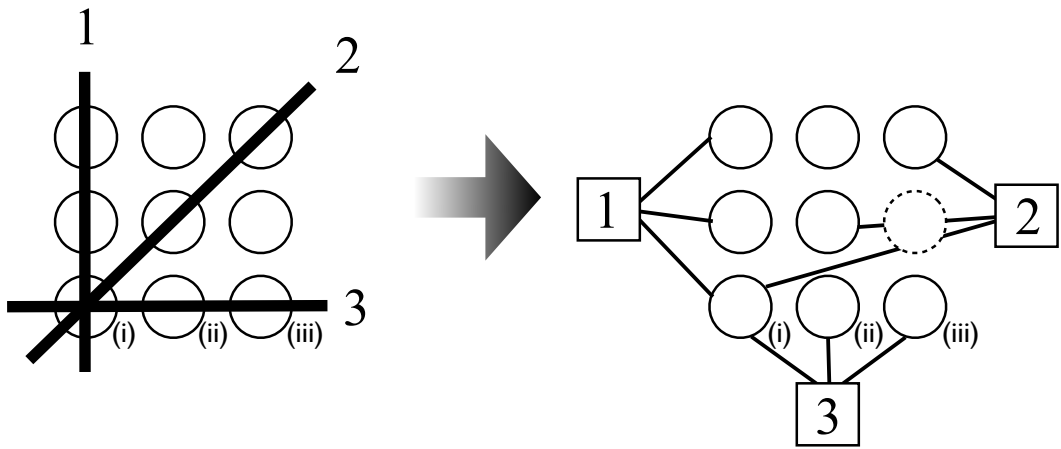}
\end{center}
\caption{Construction of the factor graph for the tomographic
  reconstruction. Each tomographic line in the image (left)
  corresponds to a factor node in the corresponding graphical model
  (right). Each of the $1,\ldots,M$ factor nodes imposes the
  statistical weight in eq.~(\ref{Bayes_Intro}). Here, for instance,
  the variables denoted (i),(ii) and (iii) are all linked to the
  factor node number $3$.\label{Howto}}
\end{figure}

Let us first reformulate the problem as a graphical model
\cite{Pearl82,MezardMontanari09}. We have $N$ variables $\s_i$ and $M$
factor nodes, each one implementing
the weight $\delta(y_{\mu} -
\sum_{i \in \mu} \s_i) e^{J_{\mu} \sum_{(ij) \in \mu} \s_i \s_j}$. The
corresponding factor graph
is shown in Fig.~\ref{Howto}.
\begin{figure}
\begin{center}
  \includegraphics[width=14cm]{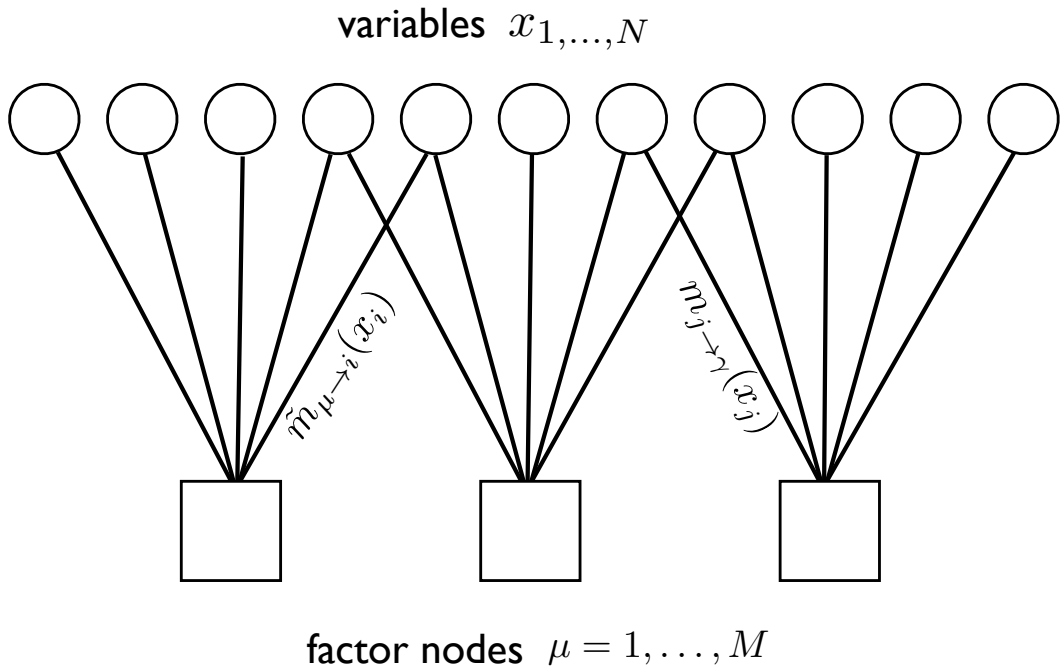}
\end{center}
\caption{Message passing on the graphical model. Each variable $j=1,\ldots,N$
  sends a message $m_{j \to \gamma} (\s_i)$ to all the
  factor nodes $\gamma$ that are connected with it, and each factor node
  $\mu=1,\ldots,M$ sends a message $\tilde m_{\mu \to i}
  (\s_i)$ to all variables $i$ that are connected
  with it.\label{GraphicalModel}}
\end{figure}

We shall not review here the derivation of the belief-propagation
algorithm, and refer instead to classical books and articles on the
subject \cite{YedidiaFreeman03,MezardMontanari09}. There are essentially
two equivalent ways to derive the BP recursion: one is to see it as a
recursion that would be exact if the factor graph were a tree (i.e. if it
had no loop). Another one is to see it as an Ansatz in the usual
variational Bayesian inference approach that usually improves on the
factorized ''mean-field'' one, leading to an expression for the
log-likelihood whose maximization leads to the belief-propagation
equations. The result is  a set of message passing equations relating
messages that go from variables to nodes and from nodes to variables. The
message $m_{i \to \gamma}$ from a variable $i$ to a node $\gamma$  is the
marginal probability of the variable $x_i$ in absence of the constraint
$\gamma$ (see Fig.~\ref{GraphicalModel}). On the other hand, the message
$\tilde m_{\mu \to i}$ from the factor $\mu$ to the variable $i$ is the
marginal probability of the variable $x_i$ when only constraint $\mu$ is
present (see Fig.~\ref{GraphicalModel}). The message $m_{i \to \gamma}$
is easily written in terms of the messages $ \tilde m_{\mu \to i}$:
\be
m_{i \to \gamma} (\s_i) \propto \prod_{\mu \in i \neq \gamma} \tilde m_{\mu
  \to i}  (\s_i)
\label{Var2Node}
\ee
where $\mu \in i$ denotes all the factor nodes (tomographic lines)
that contain variable $i$.
Note that $m_{i \to \gamma} (x_i)$
should be normalized , in the sense that $\sum_x m_{i \to \gamma} (x)=1$ (we prefer however not to write explicitly the
normalisation factor and will therefore use the sign $\propto$ instead
of $=$).  The message from a factor node to a variable, on the other hand,
reads:
\bea
\tilde m_{\mu \to i} (x_i) &\propto& \sum_{x_j; j \in \mu \neq i} \delta(y_{\mu} - \sum_{j \in \mu} \s_j)
e^{J_{\mu}    \sum_{(jk) \in \mu} \delta_{\s_j,\s_k}} \prod_{j \in \mu \neq i} m_{j
  \to \mu} (x_j) 
\label{Node2Var}
\eea
Equations (\ref{Var2Node}) and (\ref{Node2Var}) build a closed set of
equations that should be solved for the values of the
messages. Usually one seeks a solution by iterating these equations
starting from some randomly chosen initial condition.
Assuming that a fixed point of this iteration is reached, one obtains
a set of messages that solve the BP equations. The BP
estimate for the marginal distribution of $x_i$ is then given by
\be
m_i (x_i) \propto \prod_{\mu \in i} \tilde m_{\mu \to i}  (x_i)
\label{Marginal}
\ee
Let us now take a closer look at the recursion: While the iteration from
variables to factors (in eq.~(\ref{Var2Node})) is easy to compute (a
simple multiplication), the one from nodes to variables (in
eq.~(\ref{Node2Var})) is unfortunately more complicated, as it involves a
sum over all the $q$ possible values of each variable $\s$ in the line
${\mu}$. Moreover, we have to do it for each of the $N_{\mu}$ variables
in this line, and for each value it can take. At first sight, this seems
to give a complexity of order $N_{\mu} q^{N_{\mu}}$ for each iteration.
However, one can compute the message in (\ref{Node2Var}) in a much more
efficient way, namely in $O(q N_{\mu})$ iterations. Indeed, in physics
language, one can recognize that estimating $\tilde m_{\mu \to i} (x_i)$
in Eq. (\ref{Node2Var}) is nothing else that estimating the marginal of
$x_i$ in a one-dimensional graphical model known in physics as the Potts
ferromagnet in a random field at fixed total magnetization, which can be done
efficiently by transfer matrix techniques, which are equivalent for this
problem to using BP. The solution is thus to use BP in order to compute
the messages of the initial BP problem.

\subsection{Computing messages: belief propagation within belief
  propagation}

Let us now focus on the following sub-problem: we want to estimate
the marginal of a variable $\s_i$ coupled in a chain of $N_{\mu}$
variables, with a statistical weight given by
\be
P_{\mu \to i}(\vecs) \propto 
\delta(y_{\mu} -
\sum_{j \in \mu} \s_j) e^{J_{\mu} \sum_{(jk) \in \mu} \delta_{\s_j,\s_k}}
\prod_{j \in \mu \neq i} m_{j \to \mu} (x_j) \, .
\label{Pmicroc}
\ee 
There are two methods to deal with the delta function which we shall
refer to (using the statistical physics language) as the "canonical
method" and the "micro-canonical method''. The exact solution is to
use the microcanonical method where one is indeed summing only over the
configurations of the variables that satisfy the delta function in
Eq. (\ref{Pmicroc}). This can be done efficiently using a transfer
matrix approach, and we have thus implemented this method. It turns
out, however, that the canonical method is faster. It is also a
more
 natural approach when there exist errors or noise in the
measurement process. We shall therefore use the latter method in the
following.

\begin{figure}
\begin{center}
  \includegraphics[width=14cm]{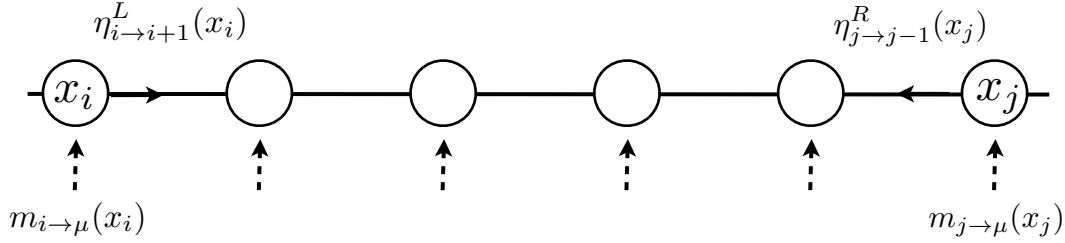}
\end{center}
\caption{Message-passing for the solving the sub-problem of estimating
  the marginal on each constraint. The left-moving message and the
  right-moving message should be computed on a one-dimensional
  line.\label{BPinBP}}
\end{figure}

The canonical method amounts to relaxing the $\delta(y_{\mu} -\sum_{j
  \in \mu} \s_j)$ constraint, and treating it only as an approximate
constraint that should hold on average, up to some small
fluctuations. More precisely, in the canonical method we replace the
delta function using a Lagrange multiplier $H$ (the "field'' in
physics language) that we shall fix later on in order to enforce the
constraint on average. The statistical weight thus becomes, instead of
(\ref{Pmicroc}):
\be
P^H_{\mu \to i}(\vecs) \propto 
e^{- H (y_{\mu} -\sum_{j \in \mu} \s_j)} e^{J_{\mu} \sum_{(jk) \in \mu} \delta_{\s_j,\s_k}}
\prod_{j \in \mu \neq i} m_{j \to \mu} (x_j)\, ,\label{cavprob}. \ee 
The value of $H$ is computed using $\mathbb{E}^H(\sum_{j \in \mu}
\s_j)=y_\mu$, where $\mathbb{E}^H$ denotes the expectation value with
respect to the measure $P^H$. Of course, the value of $H$ depends on both
$\mu$, and $i$. Now, it turns out that one can compute $H$ with an error
that becomes negligible (in the limit where the number of variables
involved in the measurement $\mu$ is large), by considering the complete
problem obtained from (\ref{cavprob}) by including all the messages $m_{j
\to \mu}$.
This complete problem is defined by:
\be P^H_{\mu}(\vecs) \propto 
e^{- H
  (y_{\mu} -\sum_{j \in \mu} \s_j)} e^{J_{\mu} \sum_{(jk) \in \mu}
  \delta_{\s_j,\s_k}} \prod_{j \in \mu} m_{j \to \mu} (\s_j) \, .\ee
We now write the nested belief propagation equations for this complete
sub-problem. Since the interaction is a two-body one, one can write it
using a version that involves only one type of message, from
variables to variables (see
\cite{KrzakalaZdeborova07,MezardMontanari09} and Fig.~\ref{BPinBP}). It is also equivalent to
the so-called transfer matrix approach in statistical mechanics
\cite{MezardMontanari09}. The recursion for the BP
messages is easily written. In order to keep notations simple, let us assume that the numbering of
variables is such that variable $i$ is neighbour to $i+1$ along direction
$\mu$. 
We shall denote the messages going from 
right to left of the line as the right-moving messages $\eta^R_{i \to
  i-1}({\s_i})$, $i=2,\dots,N_\mu$ and the ones going from  left to  right as the
left-moving messages $\eta^L_{i \to i+1}({\s_i})$, $i=1,\dots,N_\mu-1$
(see Fig.~\ref{BPinBP}). $\eta^R_{i \to
  i-1}({\s_i})$ is the marginal of $x_i$ in the absence of the left part
of the chain beyond $x_i$.
Then one has (for the left-moving messages): 
\bea
\eta^L_{i\to i+1}({\s_i}) 
&=& \frac {m_{i \to \mu} (x_i) }Z e^{H \s_i} \left[\eta^L_{i-1 \to
  i}({\s_i}) e^{J_{\mu}} + \sum_{\s_{i-1} \neq \s_i} \eta^L_{i-1 \to
  i}({\s_{i-1}}) \right] \nonumber \\ 
&=& \frac {m_{i \to \mu} (x_i) }Z e^{H \s_i} \left[1+\eta^L_{i-1\to
  i}({\s_i}) (e^{J_{\mu}}-1) \right]\, , \label{LeftMessage} \eea
where we have used that $\sum_{\s_i} \eta^L_{i\to i+1}({\s_i})=1$ and $Z$
is the corresponding normalization factor. The recursion for the
right-moving messages is similar, with the $i+1$ and $i-1$ playing
reverse roles.

Once all the $\eta^L$ and $\eta^R$ are computed, one can then obtain
the desired $\tilde m_{\mu \to i} (x_i)$, taking care of the fact that
the message $m_{i \to \mu} (x_i)$ should {\it not} be included (see
eq.~\ref{Node2Var}), using
\be \tilde m_{\mu \to i} (x_i) = 
\frac {e^{H \s_i} }Z \left[1+\eta^L_{i-1\to i}({\s_i}) (e^{J_{\mu}}-1)
\right] \left[1+\eta^R_{i+1\to i}({\s_i}) (e^{J_{\mu}}-1) \right]\, ,\ee
where $Z$ is again a normalization. 

Finally, the expected value of the variable $y_{\mu}=\sum_{i \in
  {\mu}} \s_i$ can be computed using the average value of $\s_i$ in
presence of the random field (and this time therefore including the
term $m_{i \to \mu}$) given by
\be 
\nonumber
\mathbb{E}(\s_i) = \frac{\sum_{\s_i} \s_i \tilde m_{\mu \to i} (\s_i)
m_{i \to \mu} (\s_i)}{   \sum_{\s_i}  \tilde m_{\mu \to i} (\s_i)
m_{i \to \mu} (\s_i) }\, ,
\label{field_Q}
\ee
from which we can compute the expected value of the sum of the variables $M(H)$ along
the line $\mu$:
\be
M(H)=\mathbb{E}\left(\sum_{i \in {\mu}} \s_i\right) = \sum_{i \in
  {\mu}} \mathbb{E}(\s_i) \, .
\ee
A standard result of statistical physics is that $M(H)$ is a monotonous
function of $H$. This is a consequence of the fact that there are no
phase transitions in a one-dimensional Potts model
\cite{OneDimension_66,OneDimension_Landau}.  This means that one can
apply safely a dichotomy algorithm to find the value of $H$ such that
$y_\mu=M(H)$. 

With the canonical method, one cannot fix exactly $\sum_{i \in
  {\mu}} \s_i$ but rather its expectation, and this means that we are
allowing fluctuations (typically of order $O(\sqrt{y_{\mu}})$) around
the strict tomographic constraint. Of course, one might also want to
enforce strictly this constraint, in which case one has to resort to
the microcanonical method.

Given these equations, the complete marginals for each variable can
be computed using Eq. (\ref{Marginal}). The marginal $m_i(x_i)$ denotes
the probability that variable $i$ is assigned the value $x_i$; and the
most likely labeling is thus given, at each step of the
algorithm, by
\be x_i^* = \rm{argmax}_{x_i} m_i(x_i) .
\label{Marginalization}
\ee
After convergence of the algorithm, the values $x_i^*$ thus represent the
segmentation into discrete labels that is looked for.

\subsection{A special case: binary tomography}

We shall now specialize for the rest of this paper to binary tomography,
that is, the case $q=2$, where the variables  take values $x_i=-1,1$
(after a suitable rescaling of experimental data) \footnote{Even if the
exact values of the discrete levels are not known, a method for
estimating the discrete values has been proposed by Batenburg et
al.~\cite{Batenburg2011c}}. In
physics terminology, this amounts to the Ising model, and we shall
sometimes use the term ``spins'' for pixels. In this case, it is
convenient to rewrite the expressions in a simpler way. Indeed, for
binary variables $x$ one can write a probability using a real variable
$h$ (which we will call the field) via:
\be
P(\s)=\frac{e^{h\s}}{2 \cosh{h}}\, .
\ee
Therefore, we define the fields $h_{i \to \mu}$ and $\tilde{h}_{\mu \to
i}$ by
\bea
m_{i \to \mu} (x_i)=\frac{e^{h_{i \to \mu}
    x_i}}{2\cosh{{h_{i \to \mu}}}}\, ,\\
\tilde m_{\mu \to i} (x_i)=\frac{e^{\tilde{h}_{\mu \to i}
    x_i}}{2\cosh{{\tilde{h}_{\mu \to i}}}}.
\label{def_fields}
\eea
 
Using this notation, we can rewrite the
message from a variable to a node eq.~(\ref{Var2Node}) as follows:
\be
h_{i \to \gamma} = \sum_{\mu \in i \neq \gamma} \tilde h_{\mu  \to
  i}\, .
\label{Var2Node_h}
\ee
The message from a node to a variable, on the other hand, reads:
\bea
\tilde m_{\mu \to i} (x_i) &\propto& \sum_{x_j; j \in \mu \neq i} \delta(y_{\mu} - \sum_{j \in \mu} \s_j)
e^{J_{\mu}    \sum_{(jk) \in \mu} \s_j \s_{k}} \prod_{j \in \mu \neq i} m_{j
  \to \mu} (\s_j) \\
&\propto& \sum_{x_j; j \in \mu \neq i} \delta(y_{\mu} - \sum_{j \in \mu} \s_j)
e^{J_{\mu}    \sum_{(jk) \in \mu} \s_j \s_{k} + \sum_{j \in \mu \neq i}
\s_j  h_{j \to \mu}} \, .
\eea
The marginal $\tilde m_{\mu \to i} (x_i)$ thus corresponds to the one
of a one-dimensional random field Ising model, where the external
fields are given by $h_{j \to \mu}$, and constrained to have a total
magnetization $y_{\mu}$. The recursion ``within a line'' (corresponding
to Eq.~(\ref{LeftMessage}) in the general case) can also be
written in terms of fields, using $\eta^L(\s)=\frac{e^{u^L\s}}{2
  \cosh{u^L}}$ and $\eta^R(\s)=\frac{e^{u^R\s}}{2 \cosh{u^R}}$, we have
\bea
u^L_{i\to i+1} = \atanh{\left[\tanh{J_{\mu}}\tanh{(H+h_{i\to
      \mu}+u^L_{i-1\to i})}\right]}\, ,\label{LeftMessage_g1}\\
u^R_{i\to i-1} = \atanh{\left[\tanh{J_{\mu}}\tanh{(H+h_{i\to
      \mu}+u^R_{i+1\to i})}\right]}\, .
\label{LeftMessage_g2}
\eea
Again, once these are computed, the marginal can be expressed as
\bea
\tilde m_{\mu \to i} (x_i)=\frac{e^{\tilde{h}_{\mu \to i}
    x_i}}{2\cosh{{\tilde{h}_{\mu \to i}}}}\, ,\\
{\tilde{h}_{\mu \to i}}= u^L_{i-1\to i} + u^R_{i+1\to i} + H \, .
\label{Node2Var_h}
\eea
Finally, the expected value of the variable $y_{\mu}=\sum_{i \in
  {\mu}} \s_i$, that is needed to fix the correct value of $H$, is given by summing the average value of $x_i$ (with
all random fields), each of them given by 
\be
\mathbb{E}(x_i)=\tanh{(h_{i \to \mu}+\tilde{h}_{\mu \to i})} \, ,
\label{field_Q2}
\ee
which is the analog of eq.~(\ref{field_Q}) for the binary case.
Finally, the most likely assignment for the values $x_i$ at each
step of the algorithm is given by
\be
x_i^* = \rm{sign} \((\tanh{\sum_{\mu \in i} \tilde h_{\mu \to \i}}\))
= \rm{sign} \((\sum_{\mu \in i} \tilde h_{\mu \to \i}\)).
\label{Marginalization_2}
\ee

\subsection{Implementation of the binary algorithm}
We shall now summarize the complete algorithm for the binary
case. Let us start by a few comments about the practical implementation.
\begin{itemize}
\item An efficient way to initialize the fields $\tilde{h}_{\mu \to i}$ is
  to compute their equilibrium value for $J=0$:
\be
\tilde{h}_{\mu
    \to i}=\atanh\left(y_{\mu}/N_{\mu}\right),
\ee
 where $N_{\mu}$ is the number of
  variables on the line $\mu$. This initialization is also used
in~\cite{Roux2012}.
\item Along a light ray slanted with respect to the pixel grid, the
distance between successive spins can be larger than 1. In this case, we
use a coupling between the spins that decreases when the distance between neighbors increases. We use
  $J_\mu=\atanh\left(\tanh^D{J}\right)$ where $D$
  is the $\ell_1$ distance (Manhattan or taxi-cab distance) between
  successive variables on the line $\mu$. This choice, suggested by the
  so-called high-temperature approximation \cite{Feynman}, mimics the
  effective correlation between two variables at distance $D$ in the
  two-dimensional Ising model. 
\item We have observed that the results $x_i^*$ depend only weakly on the
value of $J$, as long as $J=O(1)$. Hence, we did not optimize the value of
$J$: all the results presented here were obtained with $J=0.2$.
\item In the case of noise-free measures, we stop the algorithm when all
line-sums constraints are satisfied by the $x_i^*$. For noisy measures,
we stop the algorithm when the number of flipping spins
saturates; this criterion is discussed in Sec.~\ref{sec:noise}. 
\item We have observed empirically that damping the update of the fields 
is sometimes necessary for a proper convergence of the algorithm. We thus replaced eq.(\ref{Node2Var_h}) by
  \be {\tilde{h}_{\mu \to i}}^{t+1}= s {\tilde{h}_{\mu \to i}}^{t} +
  (1-s) \((u^L_{i-1\to i} + u^R_{i+1\to i} + H\)) \, .
\label{Node2Var_h_DAMPED}
\ee 
In our implementation, we have fixed empirically the damping value $s$
to $1-1.6/n_{\theta}$. This scaling was found to ensure convergence
for all the values of $n_\theta$ and $N$ tested.
\item In order to avoid overflow or underflow problems, we clip the
value of the different fields inside an interval $[-B, B]$. We set $B$ to
400, a value large enough for the result to be independent of $B$.  
\item In order to find the value of the external field $H$, a few
dichotomy iterations are enough to find $H$ up to a good precision
$\epsilon$ (we use $\epsilon=0.05$). In order to speed up the computation
of $H$, we use Newton's method when the total magnetization is close
enough to $y_\mu$.
\end{itemize}

\begin{codebox}
\Procname{$\proc{BP-TOMO}(y_\mu, {\rm criterion},t_{\rm max})$}
\li Initialize the messages $\tilde h$ as described above.
\li \While ${\rm conv}>{\rm criterion}$ and $t< t_{\rm max}$:
\li     \Do 
\li             $t \gets t+1$.
\li		\For each ${\mu};i \in \mu$:
\li			\Do Update $h_{i \to \mu} $ according to eq.~(\ref{Var2Node_h}).          \End
\li             \For each ${\mu};i \in \mu$:
\li                     \Do
\li                     \While $\(( \left|y_{\mu} -\sum_{i \in \mu}
  \mathbb{E}\((x_i\))\right | < \epsilon \))$
\li                          \Do Compute the messages $u^L$  and $u^R$ according to eq.~(\ref{LeftMessage_g1}-\ref{LeftMessage_g2}).
\li                        Compute the new value of $\sum_{i \in
  {\mu}}  \mathbb{E}\((x_i\))$, update $H$ (e.g. by dichotomy). \End \End
\li                         Update the $\tilde{h}_{\mu \to i} $ according to eq.~(\ref{Node2Var_h_DAMPED}).
  \End 
\li                \Do Compute the new values $x_i^*$ from
eq.(\ref{Marginalization_2}) and check for convergence.\End
\li \Return signal components $x_i^*$.
\end{codebox}
Note that this algorithm is of the embarrassingly parallel type, since each line
$\mu$ can be handled separately when computing line $7$ in the above
pseudo-code. This could be used to speed up the algorithm
by solving each line separately on different cores.  A Python code
(with some C extensions) of the current version is available on
\texttt{https://github.com/eddam/bp-for-tomo}, under a free-software
BSD license.
\newpage
\section{Numerical test in the zero-noise case \label{sec:zero_noise}}

\subsection{Generation of synthetic data}

\begin{figure}
\begin{center}
  \includegraphics[width=0.75\columnwidth]{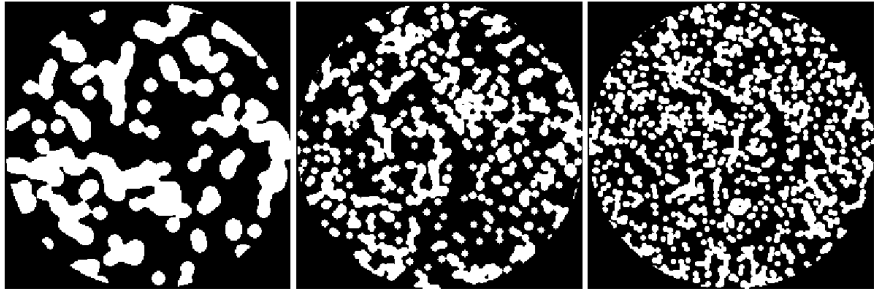}
\end{center}
\caption{Synthetic data generated for $p=14,\,28,\,38$, and a size
$256\times256$.
\label{fig:data_examples}}
\end{figure}

In order to measure the performance of our algorithm, we generate
synthetic images with different levels of sparsity of the gradients, but a 
similar binary microstructure. We obtain such a set of images with a
controlled level of gradient sparsity as follows: $p^2$ pixels are picked
at random inside the image and set to $1$ while the background is set to
$0$, then a Gaussian filter of width proportional to $1/p$ is
applied to the image. The resulting continuous-valued image is finally
thresholded at the mean intensity to obtain a binary image. Examples of
binary images with different values of $p$ are shown in
Fig.~\ref{fig:data_examples}. A simple measure of the sparsity of our images is given by
\be
\rho(x) = \frac{\ell(x)}{N},
\label{eq:density}
\ee
where $\ell(x)$ is the length of the internal boundary of one of the two
phases in the image. The internal boundary is the set of pixels that have
within their 4-connected neighbours at least one pixel belonging to the
other phase \footnote{The internal boundary can be computed by simple
mathematical morphology operations~\cite{Serra1982}, as the difference
between the origin image and its morphological erosion. Morphological
erosion is one of the two fundamental operations in mathematical
morphology~\cite{Serra1982}, realizing the intuitive idea of
eroding binary sets from their boundaries.}. $\rho$ therefore measures
the fraction of pixels lying on the boundary. We choose this measure
since the original image can be retrieved from the knowledge of the
internal boundary (and the value of a single pixel). It would also have
been possible to choose the set of pixels for which the
finite-differences gradient of the image is non-zero. However, the
latter set of pixels is redundant for retrieving the image, and leads to
a slight over-estimation of the sparsity of the image. Using
the definition of Eq. (\ref{eq:density}), we find that $\rho$ scales
linearly with $p$.

\subsection{Reconstruction in the absence of noise}

\begin{figure}
\subfigure[]{
  \includegraphics[width=0.48\columnwidth]{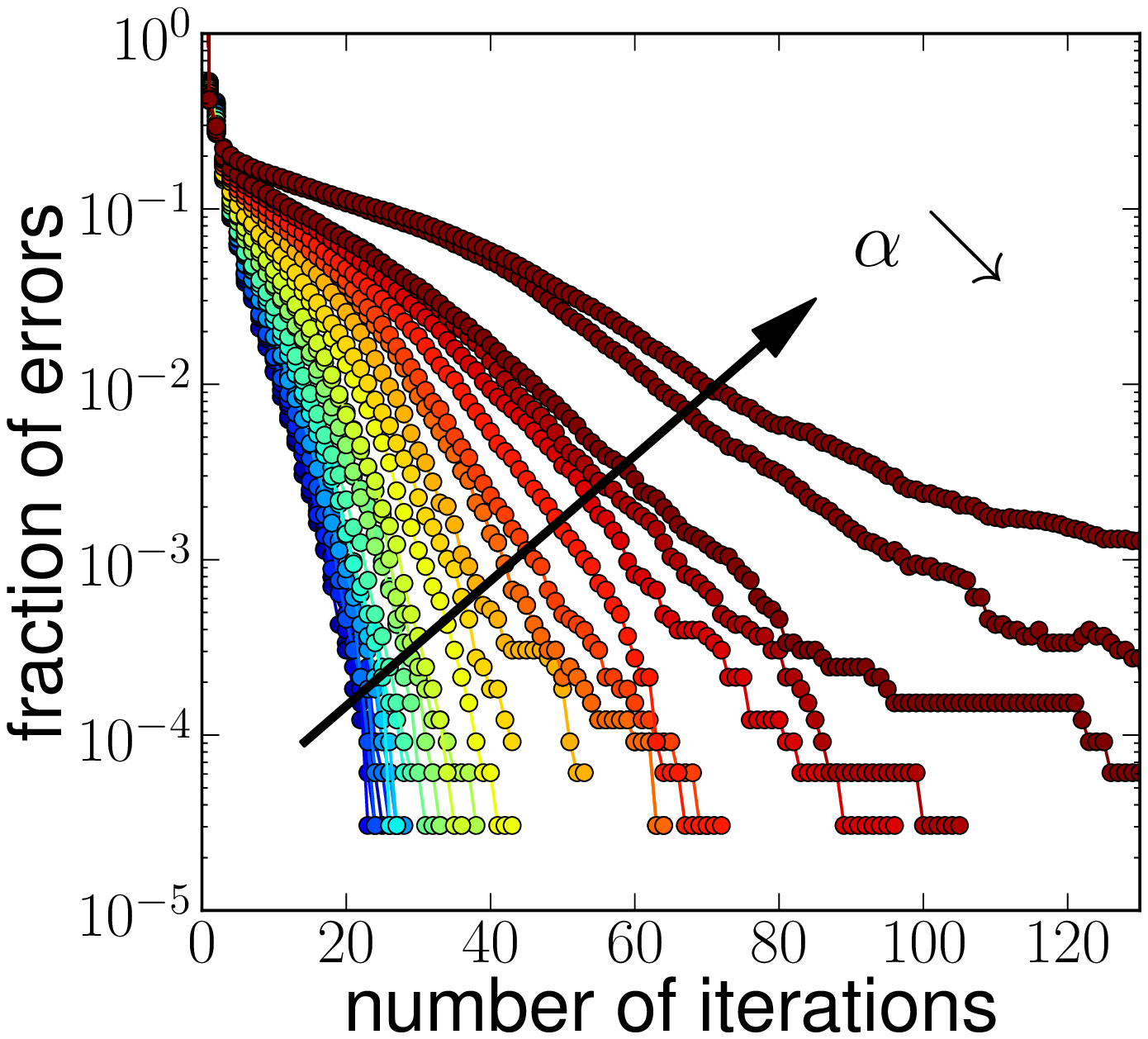}
}
\subfigure[]{
  \includegraphics[width=0.48\columnwidth]{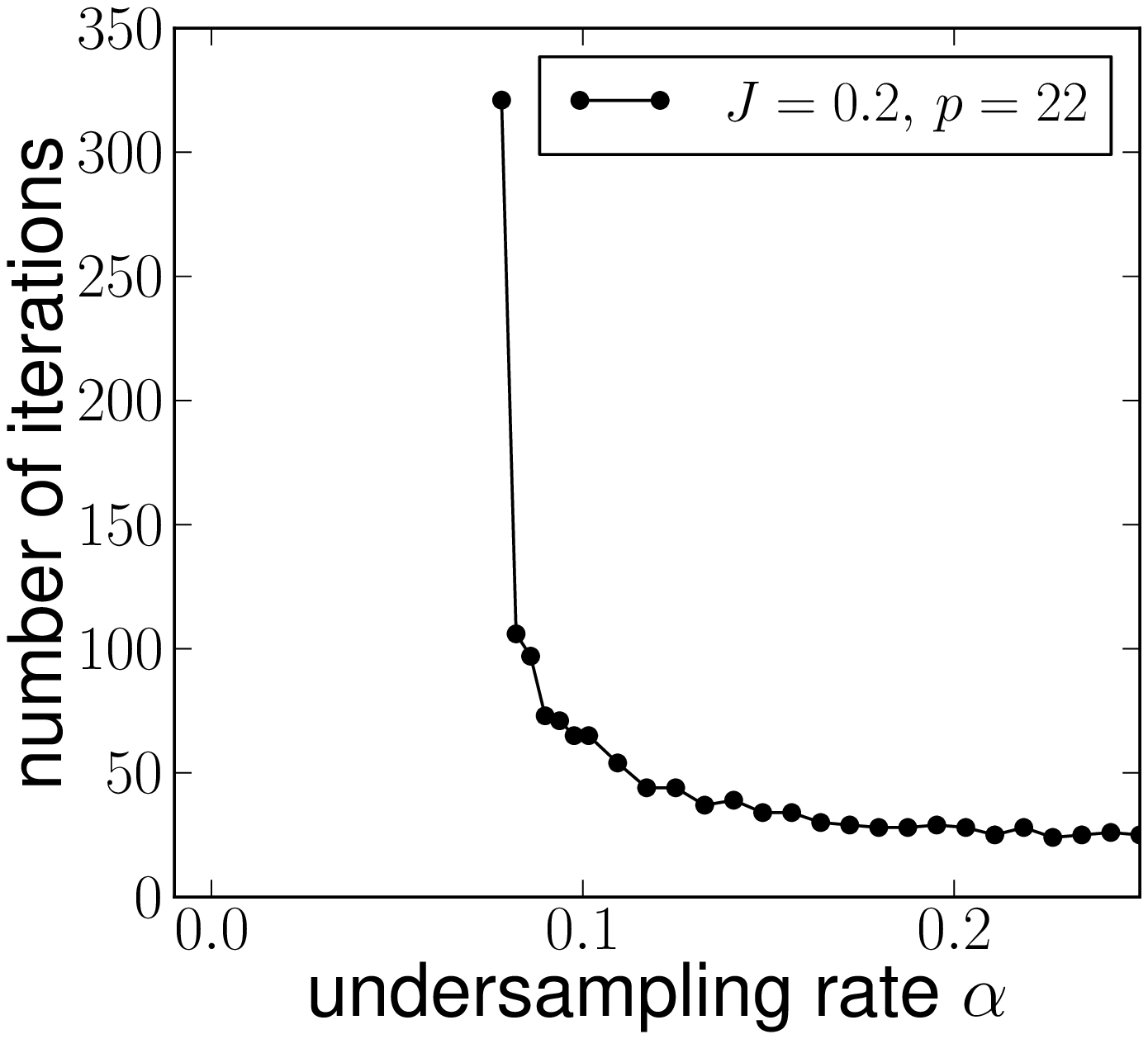}
}
\subfigure[]{
 \centerline{
\includegraphics[width=0.52\columnwidth]{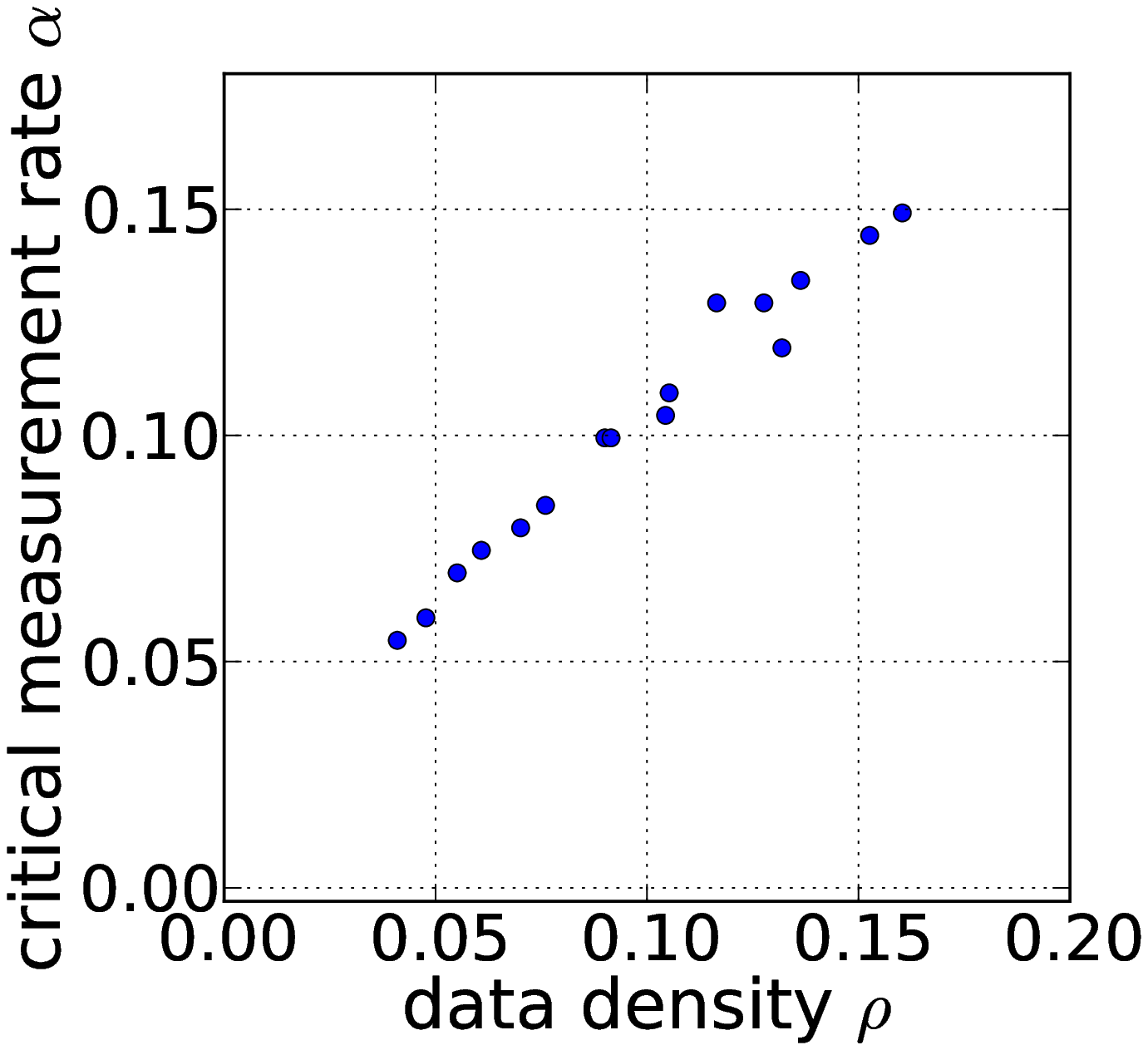}}
}
\caption{(a) Decay of the number of errors w.r.t. ground truth vs. the
number of iterations of BP-tomo, for $p=15$, $L=256$ and a number of
angles varying from 66 (blue) to 20 (red). The same value of $J=0.2$ was
used for all cases. Note how the error reduction becomes slower when
the number of projections is decreased. For smaller number of projections
(not shown here), the fraction of errors saturates at an important level.
(b) Number of iterations required to reach an exact reconstruction versus
undersampling rate $\alpha$, for $p=22$ and $J=0.2$. The number of
iterations diverges when the transition between exact reconstruction and
faulty reconstruction is approached. (c) Critical undersampling rate
$\alpha_c$ vs. boundary density $\rho$. \label{fig:phase_diagram}}
\end{figure}

Let us first investigate the effect of the number of angles on the
reconstruction of a given image, when no noise is added to the
projections: $\mathbf{y} = \mathbf{F}\mathbf{x}$.  We define the
undersampling rate by the ratio between the number of measures and the
number of variables: 
\be
\alpha = \frac{M}{N}.
\ee
Since belief
propagation is not an exact algorithm, it is not guaranteed to
reconstruct the exact original image. Nevertheless, we find
experimentally that \emph{for a sufficient number of measures}, the
algorithm always reconstructs the exact image after a finite number of
iterations. Once all spins have reached the correct orientation, we
observe that their orientation does not change any more. 
For a given
image, the evolution of
the discrepancy between the segmentation of the magnetization (Eq.
(\ref{Marginalization_2})) and the
ground truth is plotted in Fig.~\ref{fig:phase_diagram} (a) against the
number of iterations, for different measurement rates $\alpha$. The
number of errors decays roughly exponentially with time, but the decay
rate increases with the number of measurements: the convergence towards
the exact image is faster when more measurements are available. Also, a
sharp transition is observed for a critical $\alpha$, under which the
number of errors does not reach zero, and eventually increases again at
long times.  The transition between the two regimes is hard to estimate
accurately, since the convergence time seems to diverges when the
transition is approached, as shown in Fig.~\ref{fig:phase_diagram} (b).
We estimate the critical undersampling rate from the lowest number of angles
for which exact reconstruction is reached before 400 iterations.

We have measured the critical undersampling rate $\alpha_c$ for images with
different sizes of the microstructure, hence different levels of
sparsity. The critical undersampling rate is plotted against the
boundary density $\rho$ (Eq. (\ref{eq:density})) in
Fig.~\ref{fig:phase_diagram} (c). We observe a linear relationship
\be
\alpha_c \simeq \rho.
\ee
This simple relationship corresponds to a very good recovery performance,
since $\rho N$ can be seen as the number of unknowns needed to retrieve
the image. Therefore, we only need a number of measurements $M = \alpha N$
comparable to the number of unknowns for an exact reconstruction of the
image.

\newpage
\section{Numerical test in the noisy case \label{sec:noise}}

Now that we have established the phase diagram of our algorithm, we wish
to assess its performance for noisy measurements. Different sources of
noise may corrupt the measurements in X-ray tomography~\cite{Tu2006}.
Here we add a Gaussian white noise of fixed amplitude to the projections.
We define the measure noise to signal ratio (NSR) by
\be
\textrm{NSR} = \frac{\sigma}{L},
\ee
where $\sigma$ is the standard deviation of the additive noise.

\begin{figure}
\begin{center}
  \includegraphics[width=0.5\columnwidth]{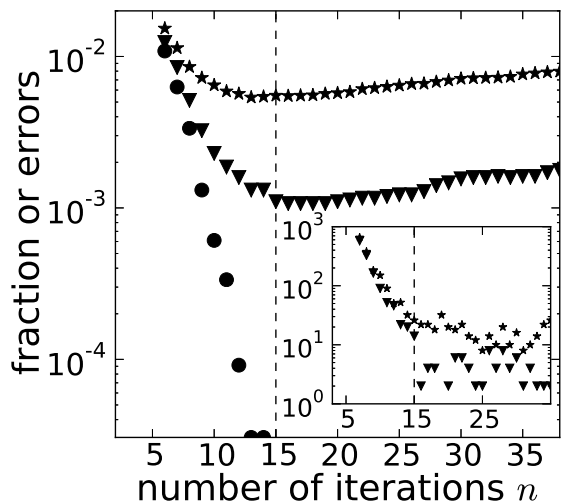}
\end{center}
\caption{\textbf{Stop criterion in the noisy case - } Main axes: for
$\alpha=1/10$, $p=14$, fraction
of error compared to ground truth vs. number of iterations $n$ for a
signal to noise ratio of 0.6 \% ($\bullet$), 1 \% ($\blacktriangledown$)
and 2 \% ($\star$). For a small amplitude of the noise, the number of
errors converges to zero ($\textrm{SNR}=0.6\%$). For a large amplitude of
the noise, the fraction of errors reaches a minimum after the number of
iterations needed to reach an exact reconstruction for a smaller value
of the noise (here, $n=15$). Inset: number of flipped spins between two
iterations vs. number of iterations. We observe that the minimum fraction
or errors correlates well with a change of slope for the number of flipped
spins.
\label{fig:stop_criterion}}
\end{figure}

In the noisy case, it is not possible any more to stop the algorithm when
all constraints are satisfied. Nevertheless, for a small value of the
noise the algorithm still converges to an exact reconstruction
(Fig.~\ref{fig:stop_criterion}). When there is convergence to the exact
solution, we observe that this solution is stable: all spins keep their
orientation during further iterations of the algorithm. Convergence is
reached after a number of iterations similar to the noise-free case. For
a larger amplitude of the noise, a finite fraction of errors remains. We
observe that the fraction of errors first decreases and reaches a
minimum, then starts to slowly increase again. A good choice of the stop
criterion is therefore important in order to optimize the quality of the
reconstruction. Empirically, we found that the best number of iterations
correlates well with the number of iterations needed to reach exact
reconstruction in the small-noise case. It can also be detected by a
change of slope in the decay of spins flipped between successive
iterations (see the inset in Fig.~\ref{fig:stop_criterion}). 

\subsection{Robustness to noise}

\begin{figure}
\begin{center}
  \includegraphics[width=0.6\columnwidth]{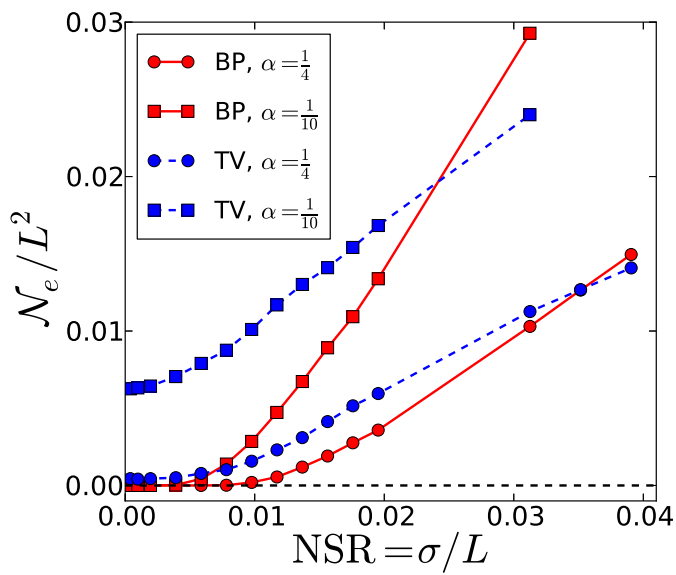}
\end{center}
\caption{Fraction of error $\mathcal{N}_e / L^2$ vs. normalized measure
noise $\sigma / L$, for BP-tomo and the convex algorithm (TV). These
numerical tests have been performed for two different undersampling rates
$\alpha$. For a small noise, BP-tomo always outperforms the convex
algorithm. In particular, there is a finite interval of noise intensity
for which the reconstruction is error-free. 
\label{fig:err_vs_sigma}}
\end{figure}

In Fig.~\ref{fig:err_vs_sigma}, we have plotted the reconstruction error
for the same image against the measure noise to signal ratio, for two
values of the measurement rate $\alpha = 1/4$ and $\alpha = 1/10$. An
image with $p=14$ (Fig.~\ref{fig:data_examples} left) was used, and
$J=0.2$ was used for all values of the NSR. We have also computed the
reconstruction error for the convex minimization of Eq.
(\ref{eq:tv_interval}). For a fair comparison, the best value of $\beta$
minimizing the reconstruction error was computed using Brent's method
(whereas no parameter was optimized for BP-tomo). We observe that the
reconstruction quality is better for BP-tomo than for the convex
algorithm, up to a threshold above which the error is greater for
BP-tomo. For the smaller undersampling rate $\alpha = 1/10$, the values
of the error reached at the threshold are too high for a satisfying
reconstruction for several applications (most pixels lying on boundaries
are wrongly labeled). Interestingly, there is a first
regime for small noise, for which the fraction of errors is zero or very
small for BP-tomo, while the error increases much faster for the convex
algorithm. For large noise, however, the convex algorithm gives better
results.

\begin{figure}
\begin{center}
\subfigure[BP-tomo]{
  \includegraphics[width=0.99\columnwidth]{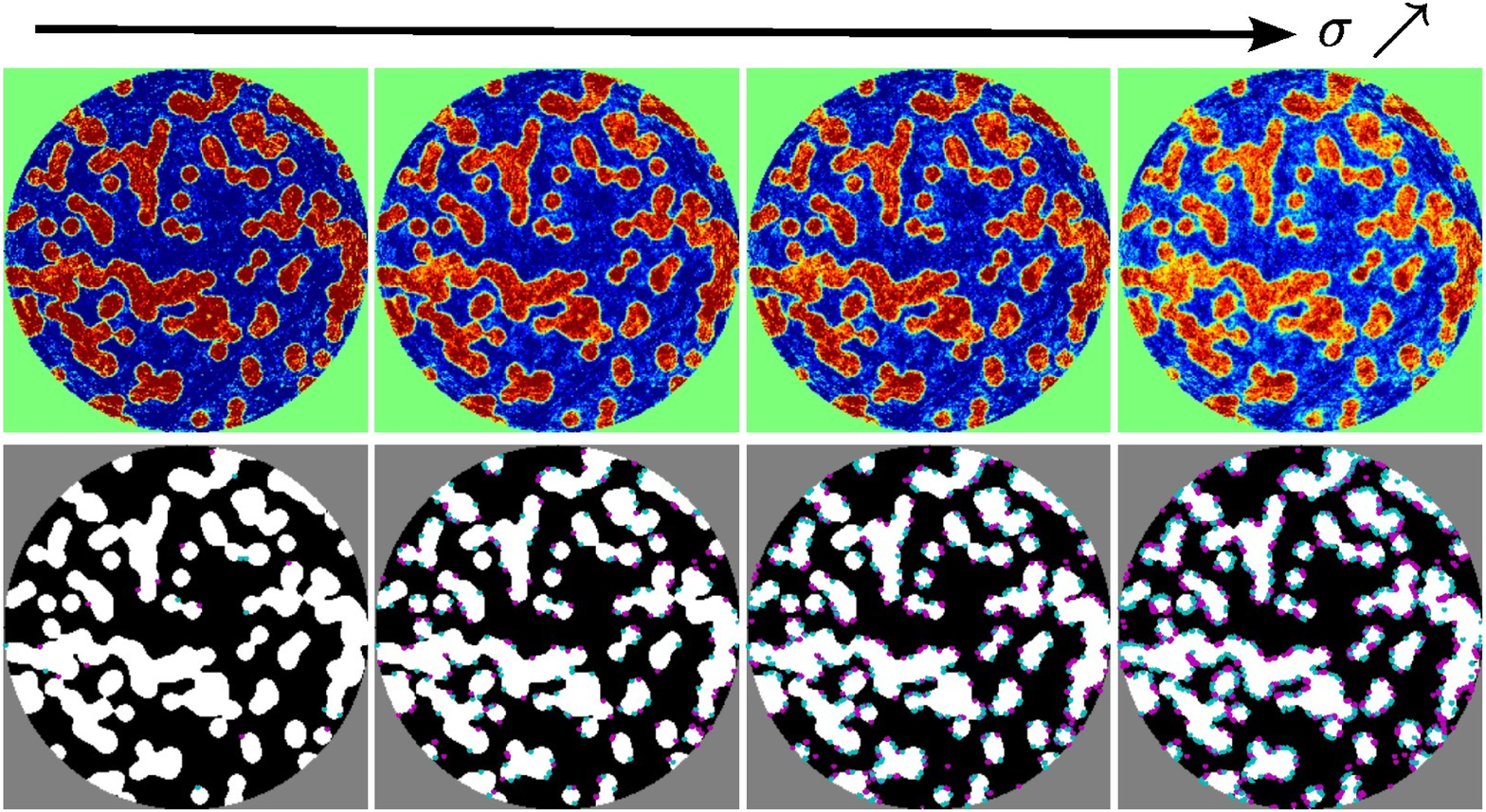}
}
\subfigure[Total Variation regularization]{
  \includegraphics[width=0.99\columnwidth]{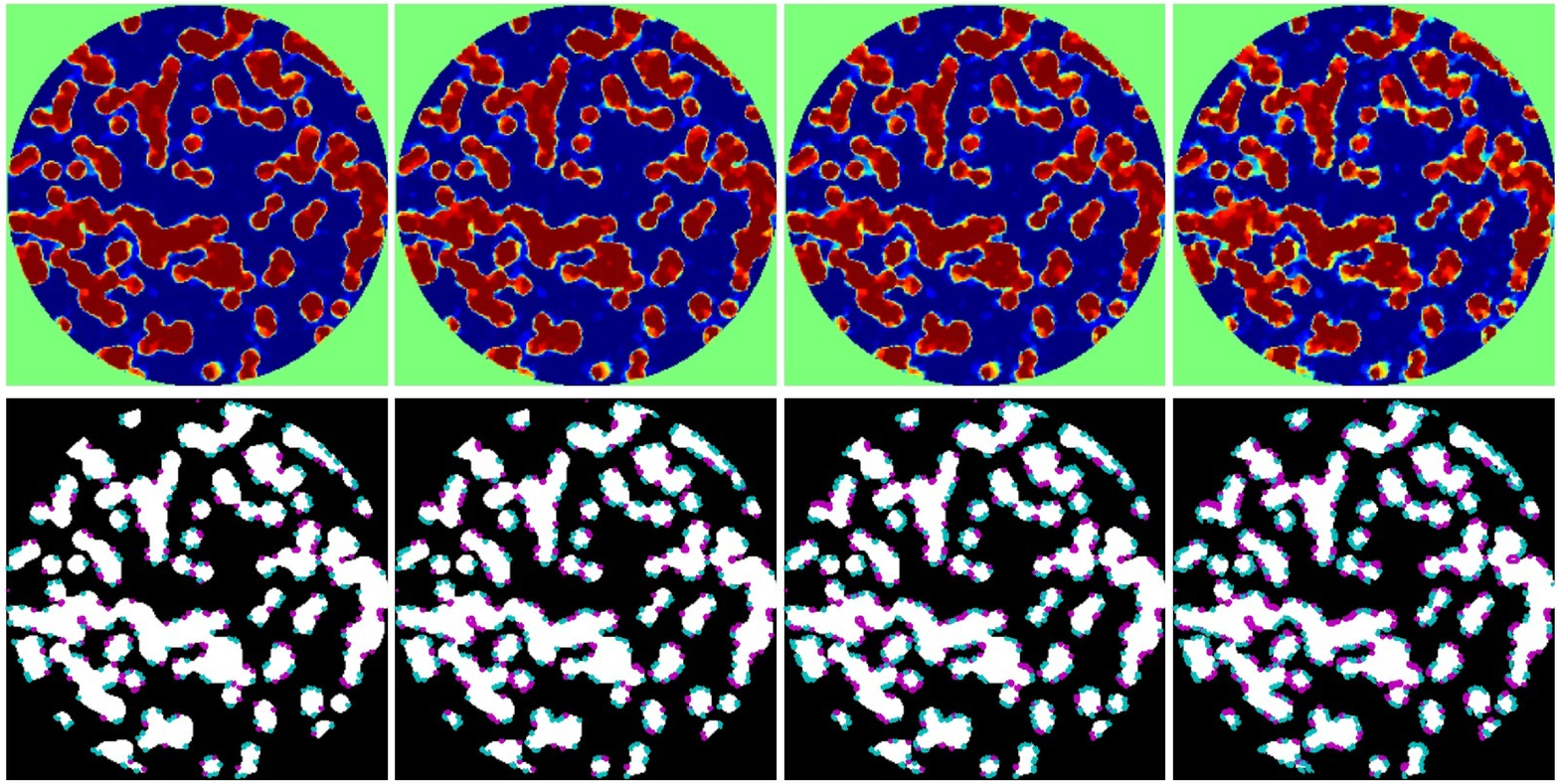}
}
\end{center}
\caption{(a) BP-tomo, $\alpha=1/10$, $p=14$, reconstructions for
$\sigma/L = 0.006,\,0.01,\, 0.02,\, 0.03$ (from left to right). The top
row shows the magnetization of the pixels, and the bottom row the
segmentation, with segmentation errors contoured in blue (resp. magenta)
for wrong pixels in the $x=-1$ (resp. $x=1$) phase. Note how the absolute
value of the magnetization decreases when the measure noise increases,
revealing a greater uncertainty. (b) TV (convex algorithm),
$\alpha=1/10$, $p=14$, reconstructions for $x/L = 0.006,\,0.01,\,0.02,\,
0.03$ (from left to right). The top row displays the result of the convex
optimization, and the bottom row its segmentation. For a low noise
amplitude, boundaries between
domains are not as sharp as for BP-tomo, resulting in more segmentation
errors close to boundaries. \label{fig:images_sigma}}
\end{figure}

Fig.~\ref{fig:images_sigma} shows that the reconstruction errors are
mostly located on the boundaries between the two phases for the two
algorithms (except for a few isolated errors for BP-tomo), but that a
larger fraction of the interfaces are correctly reconstructed for
BP-tomo. Fig.~\ref{fig:images_sigma} also displays the continuous
magnetization, as well as the non-segmented minimization of the convex
algorithm. We see that the contours of the objects are delineated more
accurately for BP-tomo than for the convex algorithm, especially for
objects with concavities.

For a given value of the NSR, we observe in Fig.~\ref{fig:err_vs_sigma}
that the reconstruction error increases when the measurement rate
decreases. The effect of the measurement rate in the noisy case is
studied in the next paragraph.

\subsection{Influence of the number of measures}

\begin{figure}
\begin{center}
  \includegraphics[width=0.6\columnwidth]{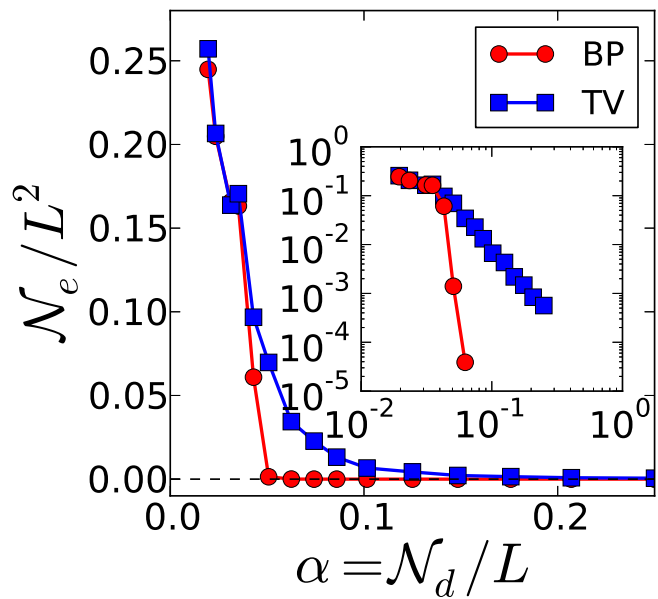}
\end{center}

\caption{Fraction or error $\mathcal{N}_e / L^2$ as a function of the
measurement rate $\alpha$, for BP-tomo and the convex algorithm (TV)
(inset in log-log coordinates). These numerical tests have been performed
for a small noise $\sigma / L = 0.002$. For both algorithms, the number
of errors decrease when the measurement rate is increased, but the
evolution is very different. For BP-tomo, the evolution is very sharp:
above a certain threshold, the reconstruction is error-free. However,
below this threshold the number of errors displays a sharp transition and
increases quickly to high values. The convex algorithm gives similar
results as BP below the BP-tomo threshold, but does not have a "phase
transition", so that it is outperformed by BP above the threshold.
\label{fig:err_vs_alpha}}

\end{figure}

For the same image ($p=14$, $J=0.2$), we have fixed the SNR to $0.002$
and computed the reconstruction error as a function as the measurement
rate, for the two algorithms. Results are shown in
Fig.~\ref{fig:err_vs_alpha}. For BP-tomo, we observe a sharp transition
between an exact reconstruction to a failed reconstruction when the
measurement rate is decreased (a score of 0.5 would be obtained when
labeling pixels at random). Therefore, the noisy case displays a
transition between failure and success as in the noise-free case
(Fig.~\ref{fig:phase_diagram}). The transition is observed at the same
measurement rate as the noise-free case. Approaching the transition from
above, a finite error fraction appears (see the log-log inset in
Fig.~\ref{fig:err_vs_alpha}). In the noise-free case, exact
reconstruction could be reached until the transition, but an increasing
number of iterations was needed when approaching the transition
(Fig.~\ref{fig:phase_diagram} (b)).

In contrast, the evolution of the error fraction when decreasing the
measurement rate is much smoother for the convex algorithm, and does not
display a phase transition in Fig.~\ref{fig:err_vs_alpha}. As a result,
the error fraction is larger for the convex algorithm above the BP-tomo
transition. At low measurement rates, the error fraction is comparable for
the two algorithms, but this regime corresponds to values of the
error fraction that are not acceptable for most applications.   

\begin{figure}
\begin{center}
  \includegraphics[width=0.6\columnwidth]{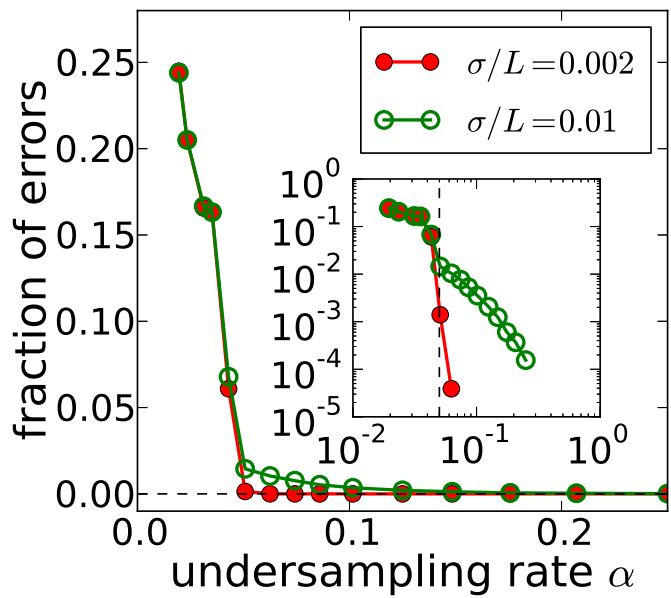}
\end{center}

\caption{Fraction or error $\mathcal{N}_e / L^2$ as a function of
the measurement rate $\alpha$, for BP-tomo and two different values of
the noise amplitude
(inset in log-log coordinates). The transition is sharper for a small
noise. In the inset, the dashed vertical line represents the empirical
value of the noise-free critical undersampling rate. 
\label{fig:alpha_sigma}}
\end{figure}

For a larger value of the noise, Fig.~\ref{fig:alpha_sigma} shows
that a transition between good and failed reconstruction is still
observed. However, the fraction of errors is larger, therefore the
transition is not as sharp as for smaller NSR. Comparing
Fig.~\ref{fig:err_vs_alpha} and Fig.~\ref{fig:alpha_sigma}, we conclude
that the reconstruction quality is better above the transition for a SNR
of $0.01$ than for the convex algorithm for a SNR five times smaller.

\newpage
\section{Discussions and perspectives}

The belief-propagation algorithm presented in this article is found to
have excellent recovery properties for tomographic binary reconstruction.
Indeed, in the noise-free case the exact original image can be
reconstructed from a number of tomographic measures approximately equal
to the number of pixels lying on the internal boundary of the objects,
that can be viewed as the number of unknowns in our problem. We are not
aware of other studies on discrete tomography where the empirical
computation of the critical undersampling rate as a function of the data
sparsity was performed. In Batenburg's DART
algorithm~\cite{Batenburg2011}, a sharp transition between exact and
faulty reconstruction is also observable for the different kinds of
phantoms studied, but the relation between the complexity of the data and
the critical undersampling rate was not described. Our search for
recovery bounds is inspired by empirical and theoretical results obtained
in the field of compressed sensing about phase transitions in
recovery~\cite{Donoho2009, Donoho2010, Krzakala2012, Krzakala2012b}.
Obtaining theoretical results for the discrete tomography problem is more
difficult than for signals sparse in a known basis -- our measure of
image ``sparsity'' comes from set theory and mathematical morphology, and
applies only to discrete images. Nevertheless, it is possible to obtain
empirical bounds as we did in Section~\ref{sec:zero_noise}. For a more
complete phase diagram of our algorithm, one should also investigate
corrections depending on the size of the data $N$, that appear for
example in the $\ell_0 - \ell_1$ Donoho-Tanner phase
transition~\cite{Donoho2009}. A probability of exact recovery could also
be obtained by computing the critical undersampling rate for a large
number of images sharing the same sparsity~\cite{Donoho2009, Donoho2010}. 

When Gaussian noise corrupts the measurements, we observe that an exact
reconstruction can still be retrieved for small noise and a sufficient
number of measurements. For a fixed non-zero level of noise, the fraction
of segmentation errors grows when the number of measurements is
decreased, but the fraction of errors is smaller than the value of the
noise over signal ratio for a large window of undersampling rates. At the
noise-free critical undersampling rate, the fraction of errors grows
very fast when the number of measurements is decreased further. Above the
noise-free critical undersampling rate and for moderate noise, we find
that our algorithm systematically outperforms a convex algorithm
implementing the convex relaxation of the binary tomography problem. The
value of the noise corresponding to the crossover when the convex
algorithm becomes better increases when the undersampling rates increases
(i.e., when one moves away from the transition). 
At the crossover, the
reconstruction error is high when the undersampling rate $\alpha$ is
small -- typically, most pixels lying on objects
boundaries are wrong --, making the reconstruction unsuitable for further
processing or interpretation in many cases. We therefore conclude that
the belief-propagation algorithm is better in most cases where further
processing of the images needs to be done. 

For practical applications, if the measurement noise is known and the
undersampling rate can be selected, we suggest that the operator should
select the undersampling rate to allow for a comfortable margin with
respect to the critical noise-free undersampling rate. This will speed
up the convergence rate (see Fig.~\ref{fig:phase_diagram} (b)) and reduce
significantly the measurement error when the noise is important (compare
the plots in Fig.~\ref{fig:err_vs_sigma} for $\alpha=1/4$ and
$\alpha=1/10$ for 1\% of noise and more). In other words, trying to approach
the limit $\alpha_c \simeq \rho$ is more interesting for theoretical than
practical reasons. However, the good news is that one does not need to go
very far from the transition $\alpha_c$ for an efficient and accurate
reconstruction, all the more for a small noise amplitude. A systematic
evaluation of the reconstruction error as a function of $\alpha,\,\rho$
and the NSR is out of the scope of this article, but would be of
interest for applications.

In future work, the recovery properties of our algorithm should be tested
on real experimental data. In our prototype code, the BP algorithm is a
few times (about 5 times) slower than the TV algorithm, principally
because of the iterations needed to find $H$. However, the number of
iterations needed to converge is much smaller for the BP algorithm, and
we expect the speed of the algorithms to depend very much on the
implementation. In order to reconstruct large 3-D volumes in a reasonable
time, our prototype implementation of the algorithm can be improved on
the numerical side -- for example with an implementation coded on the
graphics processing unit (GPU) -- and probably also on the algorithmic
side. For example, a different update scheme for the message passing
might accelerate the convergence of the algorithm. For a better precision
when the density of boundaries in the image is high (large $\rho$), the
geometrical factors $F_{\mu i}$ (see Section~\ref{sec:statphys}) should
be computed, as they are in usual algebraic reconstruction algorithms.
For 3-D images, the additional knowledge that neighboring horizontal
slices are likely to have the same pixel values can be easily implemented
in the algorithm as well. In the same way, a time-regularization in a
timeseries of several images could also be used~\cite{Pustelnik2010,
Chaari2012}, corresponding to a 3D + time regularization. This should
lead to an even better recovery compared to the 2-D case.  

We also plan to test the recovery properties in the multilabel case
($q>2$). When $q$ becomes large, it remains an open question to know how
the discrete algorithm performs compared to the convex continuous
algorithm.

Finally, for \emph{in-situ} tomography~\cite{Baruchel2006} the sample is
rotating continuously to speed-up the acquisition rate, in the case of
fast transformations. It would therefore be interesting to adapt our
algorithm to the case of continuous acquisition, when one image taken on
the detector integrates the projection of the sample over a finite
angular sector.

\section*{Acknowledgments} 
The research leading to these results was supported by the ANR program
"EDDAM" (ANR-11-BS09-027) and from the European Research Council under
the European Union's $7^{th}$ Framework Programme (FP/2007-2013)/ERC
Grant Agreement 307087-SPARCS. We gratefully acknowledge fruitful
conversations with E. Chouzenoux, S. Roux, H. Talbot and G. Varoquaux.

\newpage
\section*{References}

\bibliographystyle{unsrt.bst}
\bibliography{refs}

\begin{thebibliography}{10}

\bibitem{Slaney1988}
M.~Slaney and A.~Kak.
\newblock Principles of computerized tomographic imaging.
\newblock {\em SIAM, Philadelphia}, 1988.

\bibitem{Herman2009}
G.T. Herman.
\newblock {\em Fundamentals of Computerized Tomography: Image Reconstruction
  from Projections}.
\newblock Springer Verlag, 2009.

\bibitem{Shepp1974}
L.A. Shepp and B.F. Logan.
\newblock The {F}ourier reconstruction of a head section.
\newblock {\em IEEE Trans. Nucl. Sci}, 21(3):21--43, 1974.

\bibitem{Ramachandran1971}
GN~Ramachandran and AV~Lakshminarayanan.
\newblock Three-dimensional reconstruction from radiographs and electron
  micrographs: application of convolutions instead of {F}ourier transforms.
\newblock {\em PNAS}, 68(9):2236--2240, 1971.

\bibitem{Baruchel2006}
J.~Baruchel, J.Y. Buffiere, P.~Cloetens, M.~Di~Michiel, E.~Ferrie, W.~Ludwig,
  E.~Maire, and L.~Salvo.
\newblock {Advances in synchrotron radiation microtomography}.
\newblock {\em Scripta Materialia}, 55(1):41--46, 2006.

\bibitem{Buffiere2010}
J.Y. Buffiere, E.~Maire, J.~Adrien, J.P. Masse, and E.~Boller.
\newblock In situ experiments with {X}-ray tomography: an attractive tool for
  experimental mechanics.
\newblock {\em Experimental mechanics}, 50(3):289--305, 2010.

\bibitem{Gordon1970}
R.~Gordon, R.~Bender, and G.T. Herman.
\newblock Algebraic reconstruction techniques ({ART}) for three-dimensional
  electron microscopy and x-ray photography.
\newblock {\em Journal of theoretical Biology}, 29(3):471--481, 1970.

\bibitem{Discrete09}
G.T. Herman and A.~Kuba, editors.
\newblock {\em Discrete Tomography: Foundations, Algorithms, and Applications}.
\newblock Birkh{\"a}user, 2009.

\bibitem{Batenburg2011}
K.J. Batenburg and J.~Sijbers.
\newblock {DART}: A practical reconstruction algorithm for discrete tomography.
\newblock {\em Image Processing, IEEE Transactions on}, 20(9):2542--2553, 2011.

\bibitem{OldBayes}
B.~Carvalho, G.~Herman, S.~Matej, C.~Salzberg, and E.~Vardi.
\newblock Binary tomography for triplane cardiography.
\newblock In {\em Information Processing in Medical Imaging}, volume 1613 of
  {\em Lecture Notes in Computer Science}, pages 29--41. Springer Berlin /
  Heidelberg, 1999.

\bibitem{Sidky2008}
E.Y. Sidky and X.~Pan.
\newblock Image reconstruction in circular cone-beam computed tomography by
  constrained, total-variation minimization.
\newblock {\em Physics in medicine and biology}, 53:4777, 2008.

\bibitem{Boykov2001}
Y.~Boykov, O.~Veksler, and R.~Zabih.
\newblock Fast approximate energy minimization via graph cuts.
\newblock {\em Pattern Analysis and Machine Intelligence, IEEE Transactions
  on}, 23(11):1222--1239, 2001.

\bibitem{Gardner1995}
R.J. Gardner.
\newblock {\em Geometric tomography}, volume 2006.
\newblock Cambridge University Press Cambridge, 1995.

\bibitem{Gardner1980}
R.J. Gardner and P.~McMullen.
\newblock On {H}ammer's {X}-ray problem.
\newblock {\em Journal of the London Mathematical Society}, 2(1):171, 1980.

\bibitem{Batenburg2011b}
K.~Batenburg, W.~Fortes, L.~Hajdu, and R.~Tijdeman.
\newblock Bounds on the difference between reconstructions in binary
  tomography.
\newblock In {\em Discrete Geometry for Computer Imagery}, pages 369--380.
  Springer, 2011.

\bibitem{Aharoni1997}
R.~Aharoni, GT~Herman, and A.~Kuba.
\newblock Binary vectors partially determined by linear equation systems.
\newblock {\em Discrete Mathematics}, 171(1):1--16, 1997.

\bibitem{Liao2004}
H.Y. Liao and G.T. Herman.
\newblock Automated estimation of the parameters of gibbs priors to be used in
  binary tomography.
\newblock {\em Discrete applied mathematics}, 139(1):149--170, 2004.

\bibitem{Liao2005}
H.Y. Liao and G.T. Herman.
\newblock A coordinate ascent approach to tomographic reconstruction of label
  images from a few projections.
\newblock {\em Discrete applied mathematics}, 151(1):184--197, 2005.

\bibitem{Liao2006}
H.Y. Liao and G.T. Herman.
\newblock A method for reconstructing label images from a few projections, as
  motivated by electron microscopy.
\newblock {\em Annals of Operations Research}, 148(1):117--132, 2006.

\bibitem{Djafari2008}
A.~Mohammad-Djafari.
\newblock Gauss-markov-potts priors for images in computer tomography resulting
  to joint optimal reconstruction and segmentation.
\newblock {\em International J. of Tomography and Statistics (IJTS)},
  11:76--92, 2008.

\bibitem{Delaney1998}
A.H. Delaney and Y.~Bresler.
\newblock Globally convergent edge-preserving regularized reconstruction: an
  application to limited-angle tomography.
\newblock {\em Image Processing, IEEE Transactions on}, 7(2):204--221, 1998.

\bibitem{Weber2003}
S.~Weber, C.~Schnorr, and J.~Hornegger.
\newblock A linear programming relaxation for binary tomography with smoothness
  priors.
\newblock {\em Electronic Notes in Discrete Mathematics}, 12:243--254, 2003.

\bibitem{Schule2005}
T.~Sch{\"u}le, C.~Schn{\"o}rr, S.~Weber, and J.~Hornegger.
\newblock Discrete tomography by convex--concave regularization and {DC}
  programming.
\newblock {\em Discrete Applied Mathematics}, 151(1):229--243, 2005.

\bibitem{Weber2006}
S.~Weber, A.~Nagy, T.~Sch{\"u}le, C.~Schn{\"o}rr, and A.~Kuba.
\newblock A benchmark evaluation of large-scale optimization approaches to
  binary tomography.
\newblock In {\em Discrete Geometry for Computer Imagery}, pages 146--156.
  Springer, 2006.

\bibitem{Weber2006b}
S.~Weber, T.~Sch{\"u}le, A.~Kuba, and C.~Schn{\"o}rr.
\newblock Binary tomography with deblurring.
\newblock {\em Combinatorial Image Analysis}, pages 375--388, 2006.

\bibitem{Weber2005}
S.~Weber, T.~Sch{\"u}le, J.~Hornegger, and C.~Schn{\"o}rr.
\newblock Binary tomography by iterating linear programs from noisy
  projections.
\newblock {\em Combinatorial Image Analysis}, pages 38--51, 2005.

\bibitem{Candes2006}
E.J. Cand{\`e}s, J.~Romberg, and T.~Tao.
\newblock Robust uncertainty principles: Exact signal reconstruction from
  highly incomplete frequency information.
\newblock {\em Information Theory, IEEE Transactions on}, 52(2):489--509, 2006.

\bibitem{Candes2006b}
E.J. Candes, J.K. Romberg, and T.~Tao.
\newblock Stable signal recovery from incomplete and inaccurate measurements.
\newblock {\em Communications on pure and applied mathematics},
  59(8):1207--1223, 2006.

\bibitem{Donoho2006}
D.L. Donoho.
\newblock Compressed sensing.
\newblock {\em Information Theory, IEEE Transactions on}, 52(4):1289--1306,
  2006.

\bibitem{Boyd2004}
S.P. Boyd and L.~Vandenberghe.
\newblock {\em Convex optimization}.
\newblock Cambridge Univ Pr, 2004.

\bibitem{Combettes2011}
P.L. Combettes and J.C. Pesquet.
\newblock Proximal splitting methods in signal processing.
\newblock {\em Fixed-Point Algorithms for Inverse Problems in Science and
  Engineering}, pages 185--212, 2011.

\bibitem{Song2007}
J.~Song, Q.H. Liu, G.A. Johnson, and C.T. Badea.
\newblock {Sparseness prior based iterative image reconstruction for
  retrospectively gated cardiac micro-CT}.
\newblock {\em Medical physics}, 34(11):4476, 2007.

\bibitem{Herman2008}
GT~Herman and R.~Davidi.
\newblock Image reconstruction from a small number of projections.
\newblock {\em Inverse Problems}, 24:045011, 2008.

\bibitem{Tang2009}
J.~Tang, B.E. Nett, and G.H. Chen.
\newblock Performance comparison between total variation (tv)-based compressed
  sensing and statistical iterative reconstruction algorithms.
\newblock {\em Physics in Medicine and Biology}, 54:5781, 2009.

\bibitem{Sidky2010}
E.Y. Sidky, M.A. Anastasio, and X.~Pan.
\newblock {Image reconstruction exploiting object sparsity in boundary-enhanced
  X-ray phase-contrast tomography}.
\newblock {\em Optics Express}, 18(10):10404--10422, 2010.

\bibitem{Jia2010}
X.~Jia, Y.~Lou, R.~Li, W.Y. Song, and S.B. Jiang.
\newblock Gpu-based fast cone beam ct reconstruction from undersampled and
  noisy projection data via total variation.
\newblock {\em Medical physics}, 37:1757, 2010.

\bibitem{Anthoine2011}
S.~Anthoine, J.~Aujol, Y.~Boursier, and C.~Melot.
\newblock On the efficiency of proximal methods in cbct and pet.
\newblock In {\em Image Processing (ICIP), 2011 18th IEEE International
  Conference on}, pages 1365--1368. IEEE, 2011.

\bibitem{Pustelnik2009}
N.~Pustelnik, C.~Chaux, and J.~Pesquet.
\newblock Parallel proximal algorithm for image restoration using hybrid
  regularization.
\newblock {\em Image Processing, IEEE Transactions on}, (99):1--1, 2009.

\bibitem{Pustelnik2010}
N.~Pustelnik, C.~Chaux, J.C. Pesquet, and C.~Comtat.
\newblock Parallel algorithm and hybrid regularization for dynamic pet
  reconstruction.
\newblock In {\em Nuclear Science Symposium Conference Record (NSS/MIC), 2010
  IEEE}, pages 2423--2427. IEEE, 2010.

\bibitem{Raguet2011}
H.~Raguet, J.~Fadili, and G.~Peyr{\'e}.
\newblock Generalized forward-backward splitting.
\newblock {\em Arxiv preprint arXiv:1108.4404}, 2011.

\bibitem{Beck2009}
A.~Beck and M.~Teboulle.
\newblock Fast gradient-based algorithms for constrained total variation image
  denoising and deblurring problems.
\newblock {\em Image Processing, IEEE Transactions on}, 18(11):2419--2434,
  2009.

\bibitem{Myers2008}
GR~Myers, TE~Gureyev, DM~Paganin, and SC~Mayo.
\newblock The binary dissector: phase contrast tomography of two-and
  three-material objects from few projections.
\newblock {\em Optics express}, 16(14):10736--10749, 2008.

\bibitem{Batenburg2007}
K.J. Batenburg.
\newblock {A network flow algorithm for reconstructing binary images from
  discrete X-rays}.
\newblock {\em Journal of Mathematical Imaging and Vision}, 27(2):175--191,
  2007.

\bibitem{Roux2012}
S.~Roux, H.~Leclerc, and F.~Hild.
\newblock Tomographic reconstruction of binary fields.
\newblock In {\em Journal of Physics: Conference Series}, volume 386, page
  012014. IOP Publishing, 2012.

\bibitem{Joseph1982}
P.M. Joseph.
\newblock An improved algorithm for reprojecting rays through pixel images.
\newblock {\em Medical Imaging, IEEE Transactions on}, 1(3):192--196, 1982.

\bibitem{Dempster}
A.P. Dempster, N.M Laird, and D.B. Rubin.
\newblock Maximum likelihood from incomplete data via the em algorithm.
\newblock {\em Journal of the Royal Statistical Society}, 38:1, 1977.

\bibitem{Pearl82}
J.~Pearl.
\newblock Reverend bayes on inference engines: A distributed hierarchical
  approach.
\newblock In {\em Proceedings American Association of Artificial Intelligence
  National Conference on AI}, pages 133--136, Pittsburgh, PA, USA, 1982.

\bibitem{KschischangFrey01}
F.~R. Kschischang, B.~Frey, and H.-A. Loeliger.
\newblock Factor graphs and the sum-product algorithm.
\newblock {\em IEEE Trans. Inform. Theory}, 47(2):498--519, 2001.

\bibitem{YedidiaFreeman03}
J.S. Yedidia, W.T. Freeman, and Y.~Weiss.
\newblock Understanding belief propagation and its generalizations.
\newblock In {\em Exploring Artificial Intelligence in the New Millennium},
  pages 239--236. Morgan Kaufmann, San Francisco, CA, USA, 2003.

\bibitem{MezardMontanari09}
M.~M\'ezard and A.~Montanari.
\newblock {\em Information, Physics, and Computation}.
\newblock Oxford Press, Oxford, 2009.

\bibitem{KrzakalaZdeborova07}
L.~Zdeborov\'a and F.~Krzakala.
\newblock Phase transitions in the coloring of random graphs.
\newblock {\em Phys. Rev. E}, 76:031131, Sep 2007.

\bibitem{OneDimension_66}
Gabor~T. Herman and Attila Kuba.
\newblock {\em Mathematical Physics in One Dimension: ExactlySoluble Models of
  Interacting Particless}.
\newblock Academic Press, 1966.

\bibitem{OneDimension_Landau}
L~D Landau and E~M Lifshitz.
\newblock {\em Statistical Physics}.
\newblock Pergamon, 1959.

\bibitem{Batenburg2011c}
KJ~Batenburg, W.~van Aarle, and J.~Sijbers.
\newblock A semi-automatic algorithm for grey level estimation in tomography.
\newblock {\em Pattern Recognition Letters}, 32(9):1395--1405, 2011.

\bibitem{Feynman}
R.~P. Feynman.
\newblock {\em Statistical Mechanics: A Set of Lectures}.
\newblock Fontrier in Physics, 1972.

\bibitem{Serra1982}
J.~Serra.
\newblock {\em Image analysis and mathematical morphology}.
\newblock London.: Academic Press., 1982.

\bibitem{Tu2006}
S.J. Tu, C.C. Shaw, and L.~Chen.
\newblock Noise simulation in cone beam ct imaging with parallel computing.
\newblock {\em Physics in medicine and biology}, 51:1283, 2006.

\bibitem{Donoho2009}
D.~Donoho and J.~Tanner.
\newblock Observed universality of phase transitions in high-dimensional
  geometry, with implications for modern data analysis and signal processing.
\newblock {\em Philosophical Transactions of the Royal Society A: Mathematical,
  Physical and Engineering Sciences}, 367(1906):4273--4293, 2009.

\bibitem{Donoho2010}
D.L. Donoho and J.~Tanner.
\newblock Precise undersampling theorems.
\newblock {\em Proceedings of the IEEE}, 98(6):913--924, 2010.

\bibitem{Krzakala2012}
F.~Krzakala, M.~M\'ezard, F.~Sausset, Y.~F. Sun, and L.~Zdeborov\'a.
\newblock Statistical-physics-based reconstruction in compressed sensing.
\newblock {\em Phys. Rev. X}, 2:021005, May 2012.

\bibitem{Krzakala2012b}
F.~Krzakala, M.~M{\'e}zard, F.~Sausset, Y.~Sun, and L.~Zdeborov{\'a}.
\newblock Probabilistic reconstruction in compressed sensing: algorithms, phase
  diagrams, and threshold achieving matrices.
\newblock {\em Journal of Statistical Mechanics: Theory and Experiment},
  2012(08):P08009, 2012.

\bibitem{Chaari2012}
L.~Chaari, S.~M\'eriaux, S.~Badillo, J.-Ch. Pesquet, and P.~Ciuciu.
\newblock Multidimensional wavelet-based regularized reconstruction for
  parallel acquisition in neuroimaging.
\newblock {\em EURASIP Journal on Advances in Signal Processing}, December
  2011.
\newblock Under revision.

\end{thebibliography}

\end{document}